%
 

%
\input amstex
\magnification=\magstephalf

\documentstyle{amsppt}
\input pictex
\catcode`\@=11
\def\logo@{}
\catcode`\@=\active

\NoBlackBoxes
\TagsAsMath
\TagsOnRight
\parindent=.5truein
\parskip=3pt
\hsize=6.5truein
\vsize=9.0truein
\raggedbottom
\loadmsam
\loadbold
\loadmsbm
     
\let        \< = \langle     
\let        \> = \rangle     
\let        \| = \Vert
\define     \a          {\alpha}     
\redefine   \b          {\beta}     
\redefine   \dl         {\delta}     
\redefine   \D          {\Cal D}     
\redefine   \e          {\epsilon}     
\define     \F          {\Cal F}     
\define     \G          {\Gamma}
\redefine   \g          {\gamma}     
\define     \br         {\bold R}     
\define     \R         {\bold R}     
\define     \Z          {\bold Z}
\define     \s          {\sigma}
\define     \th         {\theta}     
\redefine   \O          {\Omega}     
\redefine   \o          {\omega}
\redefine  \w            {\omega}
\redefine   \t          {\tau}     
\define     \ub         {\subset}     
\define     \cl         {\Cal L}
\define     \E          {\Cal E}
\define     \A          {\Cal A}
\redefine     \B          {\Cal B}
\redefine   \i          {\infty}
\define     \pt         {\partial}     
\define     \bl         { L}     
\define     \ube         {\subseteq}

\define\qd{{\unskip\nobreak\hfil\penalty50\hskip2em\vadjust{}
     \nobreak\hfil$\square$\parfillskip=0pt\finalhyphendemerits=0\par}}

\topmatter

\title
The Stable Manifold Theorem for Semilinear Stochastic Evolution Equations and 
Stochastic Partial Differential Equations\\
\medskip
II: Existence of stable and unstable manifolds
 \endtitle

\author
Salah-Eldin A\. Mohammed{\footnote"$\sp {*}$"{The research of
this author is supported in part by NSF Grants DMS-9703852, DMS-9975462 and 
DMS-0203368.  \hfil\hfil}},
Tusheng\ Zhang{\footnote"$\sp {**}$"{The research of this author is supported in part 
by 
EPSRC Grant
GR/R91144.\hfil\hfil}}, and Huaizhong Zhao
{\footnote"$\sp
{***}$"{The research of this author is supported in part by EPSRC Grants
GR/R69518 and GR/R93582.\hfil\hfil \newline  
 AMS 1991 {\it subject classifications.\/} Primary 60H10, 60H20;
secondary 60H25. \newline {\it Key words and phrases.\/} Stochastic semiflow, cocycle,  
stochastic evolution equation
(see), stochastic partial differential 
equation (spde), multiplicative ergodic theorem, stationary solution, 
hyperbolicity, local stable (unstable) manifolds. 
\hfil\hfil
}}
\endauthor

\abstract
This article is a sequel to [M.Z.Z.1] aimed at completing the characterization of  the 
pathwise local structure of solutions of 
semilinear stochastic evolution equations (see's) and stochastic partial differential 
equations (spde's) near stationary
solutions. Stationary solution are viewed as random points in the infinite-dimensional 
state space, and the characterization is expressed in terms of the almost sure 
long-time 
behavior of trajectories of the equation
in relation to the stationary solution. More specifically, we establish  
{\it local stable manifold theorems\/} 
 for semilinear see's and spde's (Theorems 4.1-4.4).  
 These  results give smooth stable and unstable manifolds in the neighborhood of a 
hyperbolic stationary solution of the
underlying stochastic equation.  The stable and unstable manifolds are stationary, 
live 
in a stationary tubular neighborhood of the stationary solution  and are 
asymptotically 
invariant 
under the stochastic semiflow of the see/spde. The proof uses infinite-dimensional 
multiplicative ergodic theory techniques and interpolation arguments (Theorem 2.1). 
 \endabstract

\endtopmatter

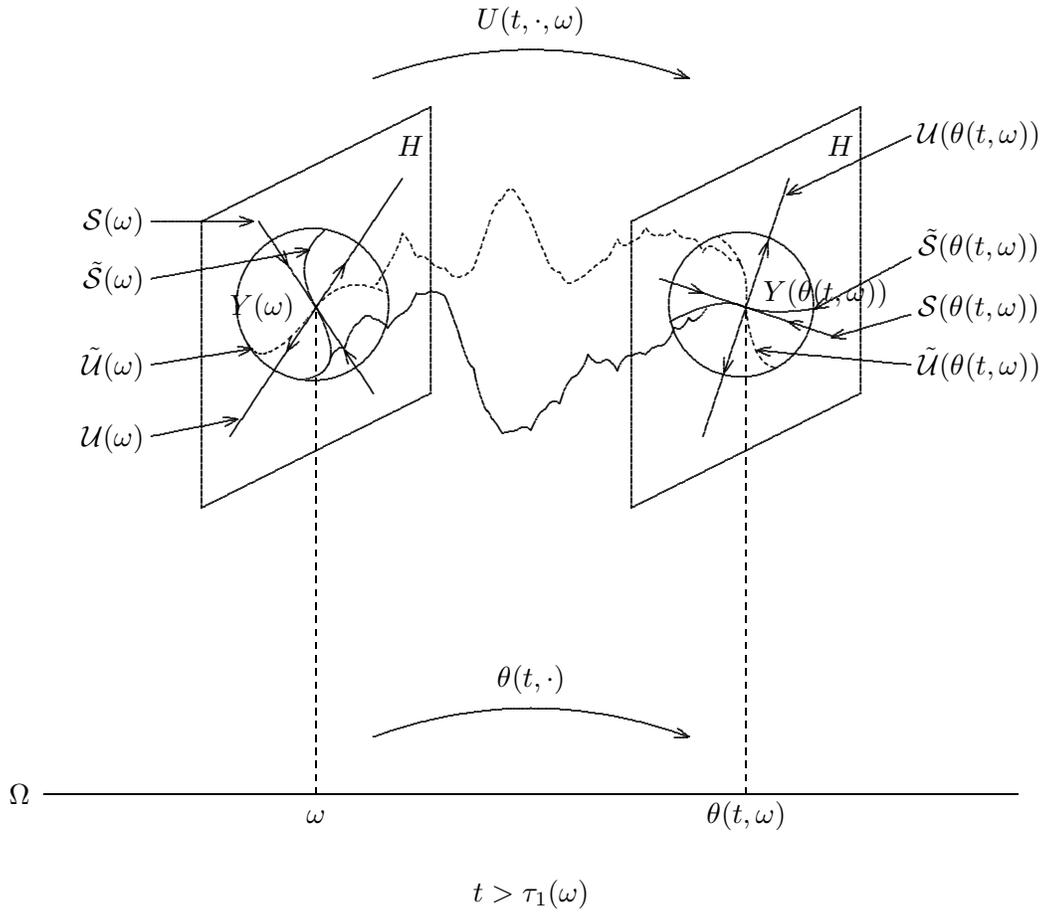
\rightheadtext{The Stable Manifold Theorem  for SEE's and SPDE's}

\leftheadtext{S.-E.A. Mohammed, T.S. Zhang and H.Z. Zhao}

\document
 
\baselineskip=16truept
 

\subheading{1. Introduction. Hyperbolicity of a stationary trajectory}

In [M-Z-Z.1], we established the existence of perfect differentiable cocycles
generated by mild solutions of a large class of semilinear stochastic evolution 
equations (see's) and stochastic partial 
differential equations (spde's). The present article is a continuation of the 
analysis in [M-Z-Z.1].  In this paper we introduce the concept 
of a stationary trajectory for the see. Within the context of stochastic differential 
equations (with memory) (sde's and
sfde's), this concept has been used 
extensively in previous work of one of the authors with M. Scheutzow ([M-S.1],
[M-S.2-4]). 
Our main 
objective is to characterize the pathwise local structure of solutions of 
semilinear see's and spde's near stationary solutions.  We introduce the concept of 
{\it hyperbolicity} for 
a stationary solution of an see. Hyperbolicity is defined by the non-vanishing of  the 
Lyapunov spectrum of the linearized 
cocycle. The hyperbolic structure of the stochastic semiflow leads to 
 {\it local stable manifold theorems\/} (Theorems 4.1-4.4)
 for semilinear see's and spde's. For a hyperbolic stationary solution of the  
see, this gives smooth stable and unstable manifolds in a neighborhood of the 
stationary 
solution. 
 The stable and unstable manifolds are stationary, live in a stationary 
tubular neighborhood of the stationary solution  and are asymptotically 
invariant under the stochastic semiflow. The proof of the stable manifold 
theorem uses infinite-dimensional multiplicative ergodic theory techniques 
([Ru.1], [Ru.2]) together with  interpolation and perfection arguments 
([Mo.1], [M-S.4]). In particular, we will assume that the reader is familiar with 
the results and the techniques in Ruelle's articles [Ru.1] and [Ru.2]. Our results 
cover 
semilinear stochastic evolution 
equations, stochastic parabolic equations, stochastic reaction-diffusion equations, 
and 
Burgers equation with additive
infinite-dimensional noise.

We recall below the definition of a cocycle in Hilbert space.
 
Let $(\O, \F, P)$  be a complete probability space. Suppose 
$\theta :\R \times \O \to \O$ is a group of $P$-preserving ergodic
 transformations on  $(\O,\F,P)$. Denote by $\bar \F$ the $P$-completion of $\F$. 

Let $H$ be a real separable Hilbert space with  norm $|\cdot|$ and   
Borel $\sigma$-algebra $\Cal B(H)$. 

Take $k$ to be any  non-negative 
integer and  $\e \in (0,1]$. Recall that a $C^{k, \e}$ {\it perfect 
cocycle\/} $(U, \theta)$  on $H$
 is a $(\Cal B (\R^+) \otimes \Cal B(H) \otimes \F, \Cal B(H))$-
measurable random field
$U:\R^+ \times H \times \O \to H$
with the following properties:
 
\item{(i)} For each $\o\in \O$, the map 
$\R^+ \times H \ni (t,x) \mapsto U(t,x,\o) \in H$
 is continuous; for fixed $(t,\o)\in \R^+ \times \O$, the map 
$H \ni x \mapsto U(t,x,\o) \in H$ is $C^{k,\e}$ ($D^k U(t,x,\o)$ is 
$C^{\e}$ in $x$ on bounded subsets of $H$).

\item{(ii)} $U(t_1+t_2,\cdot,\o)=  
U(t_2,\cdot,\theta (t_1,\o))\circ U(t_1,\cdot,\o)$
for all $t_1,t_2 \in \R^+$, all $\o \in \O$.

\item{(iii)} $U(0,x,\o)=x$ for all $x \in H, \o \in \O$.

\medskip

 We now introduce the concept of a  
{\it stationary point\/} for a cocycle $(U,\theta)$. Stationary points play the role 
of 
stochastic equilibria for 
the stochastic dynamical system. 
  %
%
 
\medskip


\definition{Definition 1.1} 

An $\Cal F$-measurable random variable
$Y:\O \to H$ is said be a {\it stationary random point} for the cocycle $(U,\theta)$ 
if  
it satisfies
the following identity:
$$ 
U(t,Y(\o),\o)= Y(\theta (t,\o)) \tag1.1
$$
for all $(t,\omega )\in \R^+ \times \O$.
\enddefinition

The reader may note that the above definition is an
infinite-dimensional analogue of a corresponding concept of
invariance that was used by one of the authors in joint work with M.
Scheutzow to give a proof of the stable manifold theorem for
stochastic ordinary differential equations (Definition 3.1, [M-S.3]).
Definition 1.1 essentially gives a useful realization of the idea
of an invariant measure for a stochastic dynamical system generated by an   
spde or a see. Such a realization allows us to
analyze the local {\it almost sure} stability properties of the
stochastic semiflow in the neighborhood
of the stationary point. The existence (and uniqueness/ergodicity) of a stationary 
random point for various classes of spde's and see's has been studied by many 
researchers.  
In this article, we will move beyond the issue of existence of stationary solutions,
and apply our stable/unstable manifold theorem to 
examine further the 
almost sure asymptotic structure of the stochastic flow generated by several 
well-known classes of see's and spde's. In particular, we establish the existence 
of local stable and unstable manifolds near their stationary points.


The main objective of this section is to define the concept of {\it hyperbolicity} 
for a stationary point $Y$ of the cocycle  $(U, \theta)$.

\medskip

First, we linearize the $C^{k,\e}$ cocycle $(U,\theta)$ along a stationary 
random point $Y$. By taking Fr\'echet derivatives at $Y(\o)$ on each side of the 
cocycle identity (ii) above,  using the chain rule and the definition of 
$Y$, we immediately see 
that $(DU(t,Y(\o),\o), \theta (t,\o))$ is an 
$L(H)$-valued perfect cocycle. Secondly, we appeal to the following classical result 
which goes back to Oseledec
in the finite-dimensional case, and to D. Ruelle in infinite dimensions.

\noindent
\proclaim{Theorem 1.1} (Oseledec-Ruelle) 
 
Let $T: \R^+ \times \O \to L(H)$ be strongly measurable, such that $(T, \theta)$ is 
an $L(H)$-valued cocycle, with each $T(t,\o)$ compact.  Suppose that  
$$ 
\displaystyle E \sup_{0\leq t \leq 1} \log^+ \|T(t,\cdot)\|_{L(H)} + 
\displaystyle E \sup_{0\leq t \leq 1} \log^+ \|T(1-t,\theta (t,\cdot))\|_{L(H)} 
< \infty.
$$  
 Then there is a sure event $\O_0 \in \F$  such 
that  $\theta (t,\cdot)(\O_0)\subseteq \O_0$ for all $t \in \R^+$, and for 
 each $\o \in \O_0$, the limit  
$$ \Lambda (\o) := \lim_{t \to \infty} [T(t,\o)^* \circ T(t,\o)]^{1/(2t)}  $$
exists in the uniform operator norm. Each linear operator $\Lambda (\o)$ is compact, 
non-negative and self-adjoint with  a
discrete  spectrum
$$ e^{\lambda_1} >  e^{\lambda_2 }  > e^{\lambda_3 }> \cdots $$
where the $\lambda_i$'s are distinct and non-random.  Each eigenvalue $e^{\lambda_i} > 
0$ has a fixed  finite
non-random multiplicity $m_i$ and  a corresponding eigen-space $F_i (\o)$, with 
$m_i := \hbox{dim} \, F_i (\o)$. Set $i=\infty$ when $\lambda_i=-\infty$. Define 
$$  E_1 (\o):= H,\quad E_i (\o):= \bigl [\displaystyle \oplus_{j=1}^{i-1}
F_j (\o) \bigr]^{\perp},\,\,  i > 1,\,\, E_{\infty}:= \hbox{ker}\, \Lambda (\o).$$
Then
$$ 
E_{\infty}\subset \cdots \subset \cdots \subset E_{i+1}(\o)\subset E_{i} (\o) \cdots 
\subset E_2 (\o) \subset E_1 (\o) =H,
$$
$$ 
\lim_{t \to \infty} \frac{1}{t} \log |T(t,\o) x| =\cases   \lambda_i  \quad
\hbox{if} \quad x \in E_i (\o) \backslash E_{i+1} (\o),\\
-\infty \quad \hbox{if} \quad x \in E_{\infty}(\o),
\endcases
$$
and 
%
$$
     T(t,\omega) (E_i (\omega)) \subseteq E_i (\theta (t, \o))
$$
for all $\, t \geq 0,\,\, i \geq 1$.

\endproclaim

\bigskip

The following figure   illustrates the Oseledec-Ruelle theorem.

\bigskip

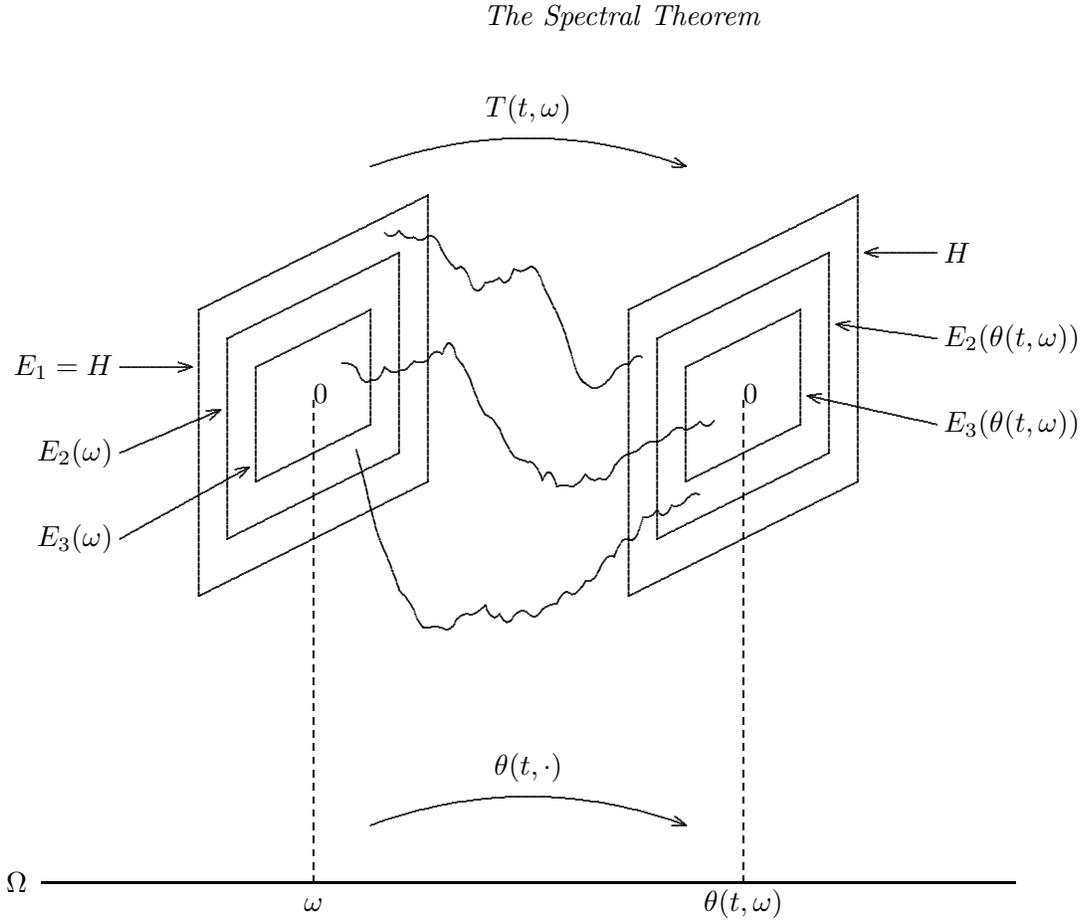
\centerline{\it The Spectral Theorem }

\bigskip
\bigskip


\beginpicture
\setcoordinatesystem units <0.075truein, 0.075truein>
\setplotarea x from -5 to 65, y from -10 to 30
%
%
\plot  8 20   24 28   24 48    8 40    8 20 /
\plot 10 24   22 30   22 44   10 38   10 24 /
\plot 12 28   20 32   20 40   12 36   12 28 /
%
%
\plot 38 20   54 28   54 48   38 40   38 20 /
\plot 40 24   52 30   52 44   40 38   40 24 /
\plot 42 28   50 32   50 40   42 36   42 28 /
%
%
\putrule from -3 0 to  65 0
\setdashes <2.5truept>
\putrule from 16  0  to  16 34
\putrule from 46  0  to  46 34
\put {$0$}                 [cb] at 16.5 33.5
\put {$0$}                 [cb] at 46.5 33.5
\put {$\omega$}                 [cb] at 16 -2.0
\put {$\theta(t,\omega)$}     [cb] at 46 -2.5
\put {$\Omega$}                 [rc] at -4 0
\setsolid
\setquadratic    
\plot 20 50   31 52   42 50 /
\arrow <.07truein> [.4, .7] from 41.9 50.034 to 42 50
\put {$T(t,\omega)$} [cb] at 31 53
\plot 20 4    31 6    42 4  /
\arrow <.07truein> [.4, .7] from 41.9 4.034 to 42 4
\put {$\theta(t,\cdot)$} [cb] at 31 7
\put {$E_1 = H$}              [rc] at  2 36
\put {$E_2(\omega)$}            [rc] at  2 30
\put {$E_3(\omega)$}            [rc] at  2 24
\put {$H$}                    [lc] at 60 44
\put {$E_2(\theta(t,\omega))$} [lc] at 60 38
\put {$E_3(\theta(t,\omega))$} [lc] at 60 32
\arrow <.07truein> [.4, .7] from        2.5 36 to 7.5 36
\arrow <.07truein> [.4, .7] from        2.5 30 to 9.5 33
\arrow <.07truein> [.4, .7] from        2.5 24 to 11.5 29
\arrow <.07truein> [.4, .7] from        59.5 44 to 54.5 44
\arrow <.07truein> [.4, .7] from        59.5 38 to 52.5 39
\arrow <.07truein> [.4, .7] from        59.5 32 to 50.5 34
%
%
\setquadratic
\plot
21 45.361
21.5 45.1653606
22.0 45.5465769
22.5 45.0465048
23.0 45.0850000
23.5 44.8808125
24.0 45.1492691
24.5 45.0225913
25.0 44.1090001
25.5 43.3528271
26.0 42.8088465
26.5 42.8839426
27.0 41.8050000
27.5 41.3461917
28.0 41.9328451
28.5 41.9935758
29.0 42.2420000
29.5 41.9657812
30.0 42.8727736
30.5 42.6928792
31.0 42.9300001
31.5 42.8570587
32.0 41.8360607
32.5 40.7600324
33.0 39.7900000
33.5 38.5362343
34.0 36.9239837
34.5 35.4477414
35.0 34.8780000
35.5 34.5560042
36.0 34.6400045
36.5 35.0030024
37.0 35.6670000
37.5 36.1129986
38.0 36.1739986
38.5 36.7099991
39.0 36.6510000
/
%
%
\plot
18 36.264
18.5 36.2305516
19.0 35.9514826
19.5 35.0451723
20.0 35.1030000  
20.5 35.2025908
21.0 36.0195523
21.5 35.6627376
22.0 35.7727000  
22.5 35.8829602
23.0 36.4153083
23.5 36.3465023
24.0 36.4900000
24.5 37.0525685
25.0 36.8622142
25.5 37.6672527
26.0 36.6690000
26.5 36.0154531
27.0 34.9983348
27.5 33.6420491
28.0 32.6400000
28.5 32.3648064
29.0 31.6629465
29.5 30.5751134
30.0 29.7590000
30.5 29.5326338
31.0 29.3523792
31.5 29.1089351
32.0 28.0710001
32.5 28.6500961
33.0 27.6210368
33.5 27.7564593
34.0 27.6180000
34.5 27.9094822
35.0 28.5164735
35.5 28.1247281
36.0 28.2569999
36.5 28.2196629
37.0 29.0505692
37.5 29.6381908
38.0 29.9440000
38.5 30.1106162
39.0 30.5022497
39.5 30.7822584
40.0 30.7920000
40.5 31.0401848
41.0 31.4969319
41.5 31.7587131
42.0 31.5710000
42.5 31.8257698
43.0 32.1980227
43.5 31.9032641
44.0 32.2980000
/
%
%
\plot
19 30.215
19.5 28.2731988
20.0 26.1711182
20.5 24.3444784
21.0 23.0880000
21.5 21.8425579
22.0 20.3966454
22.5 19.8589101
23.0 19.0930000
23.5 18.0090695
24.0 17.7863003
24.5 17.8393810
25.0 17.7980001
25.5 17.6989140
26.0 18.3581531
26.5 18.7128156
27.0 18.6819999
27.5 18.9291496
28.0 19.3910872
28.5 18.5659815
29.0 18.7950000
29.5 18.1813006
30.0 18.6389979
30.5 18.9461964
31.0 18.8120000
31.5 18.5953982
32.0 19.0899210
32.5 19.3519832
33.0 19.4850000
33.5 20.2233565
34.0 20.1363180
34.5 20.3541206
35.0 21.0930000
35.5 21.1689261
36.0 21.9408070
36.5 22.6802841
37.0 22.7420000
37.5 23.7648763
38.0 23.9739542  
38.5 24.1405551
39.0 25.4830000
39.5 25.5670688
40.0 26.0213762
40.5 25.6909954
41.0 26.3400000
41.5 26.3310987
42.0 26.7005413
42.5 27.1172134
43.0 27.0559999
/
\endpicture
%

\medskip

\noindent
\demo{Proof of Theorem 1.1}


 The proof is based on a discrete version of Oseledec's multiplicative ergodic theorem 
and the perfect 
ergodic theorem ([Ru.1], I.H.E.S Publications, 1979, pp.
303-304; cf. [O], 
[Mo.1], Lemma 5. See also Lemma 3.1 (ii) of this article). Details of the 
extension to continuous time are given in   
[Mo.1] within the context of linear stochastic functional differential 
equations. The arguments in [Mo.1] extend directly to general linear cocycles 
in Hilbert space. Cf. [F-S].
 \qd
\enddemo




\noindent
\definition{Definition 1.2 }

The sequence 
$\{\cdots < \lambda_{i+1} < \lambda_i < \cdots < \lambda_2 < \lambda_1 \}$
in the Oseledec-Ruelle theorem (Theorem 1.1) is called the 
{\it Lyapunov spectrum} of the linear cocycle $(T,\theta)$.
\enddefinition

Hyperbolicity of a stationary point $Y: \O \to H$ of the non-linear cocycle 
$(U, \theta)$ may now be defined in terms of a spectral gap in the Lyapunov 
spectrum of the linearized cocycle $(D U(t, Y(\o),\o), \theta (t,\o))$.

\noindent
\definition{ Definition 1.3}

 Let $(U,\theta)$ be a $C^{k,\e}$ ($k \geq 1, \e \in (0,1]$) perfect cocycle on a 
separable Hilbert space 
$H$ such that $U(t,\cdot,\o): H \to H$ takes bounded sets into relatively 
compact sets for each $(t, \o) \in (0,\infty) \times \O$. A stationary point $Y(\o)$ 
of the  cocycle $(U,\theta)$ is {\it hyperbolic\/} if 
\item{(a)}For any $a \in (0,\infty)$,
$$ \int_{\O} \log^+ \sup_{0\leq t_1,t_2 \leq a} \|DU(t_2, Y(\theta
(t_1,\o)),\theta (t_1,\o))\|_{L(H)} \, dP(\o) < \infty. 
$$
\item{(b)} The linearized cocycle $(D U(t, Y(\o),\o), \theta (t,\o))$
 has a non-vanishing Lyapunov spectrum 
$\{\cdots < \lambda_{i+1} < \lambda_i < \cdots < \lambda_2 < \lambda_1 \}$,
viz. $\lambda_i \neq 0$ for all $i \geq 1$.
\enddefinition

\bigskip


\medskip

By the Oseledec theorem (Theorem 1.1), the integrability condition in Definition 1.2 
(a) 
implies the existence of a discrete
Lyapunov spectrum for the linearized cocycle $(D U(t, Y(\o),\o), \theta (t,\o))$ in 
Definition 1.2 (b) above.
 
The following result is a random version of the saddle point property for 
hyperbolic linear cocycles. A proof is given in ([Mo.1], Theorem 4, Corollary 2; 
[M-S.1], Theorem 5.3) within the context of
stochastic differential systems with memory; but the arguments therein extend 
immediately to linear cocycles 
in Hilbert space. 

\proclaim{Theorem 1.2} (Stable and unstable subspaces)

Let $(T,\theta)$ be a linear cocycle on a  Hilbert space $H$. 
Assume that $T(t, \o): H \to H$ is a compact linear operator for 
each $t > 0$ and a.a. $\o \in \O$. Suppose that 
$$
 E \log^+ \sup_{0\leq t_1,t_2 \leq 1} \|T(t_2,\theta (t_1, \cdot))\|_{L(H)} < \infty,
$$
and let the cocycle $(T, \theta)$ have  a non-vanishing Lyapunov spectrum
$\{\cdots < \lambda_{i+1} < \lambda_i < \cdots < \lambda_2 < \lambda_1 \}$.
Pick  $i_0 > 1$ such that  $\lambda_{i_0} < 0 < \lambda_{i_0-1}$. 
 
Then there is a sure event $\O^* \in \F$ and {\it stable} and {\it unstable} subspaces
$\{ S(\o), U(\o): \o \in \O^* \}$, $\F$-measurable (into the Grassmanian),   such that 
for each  $\o \in \O^*$, the following is
true:

\item{(i)} $\theta (t,\cdot)(\O^*)= \O^*$ for all $t \in \R$.
\item{(ii)} $ H = {\Cal U} (\omega) \oplus {\Cal S} (\omega).$ The subspace 
$\Cal U(\o)$ 
is  finite-dimensional with a fixed non-random dimension, and $\Cal S(\o)$ is 
closed with a finite non-random codimension.
\item{(iii)}({\it Invariance})
     $$ T(t, \o) ({\Cal U} (\omega))= {\Cal U} (\theta(t,\omega)), \,\,
     T(t, \o) ({\Cal S} (\omega)) \subseteq {\Cal S}(\theta (t, \omega)),  
 $$
for all $t \geq 0$, 

\item{(iv)}({\it Exponential dichotomies}) 
$$
     |T (t, \omega)(x) |  \geq | x |
 e^{\delta _1 t} \quad  \hbox{for all} \quad \,\,t \geq \tau^* _1,
 x \in {\Cal U} (\omega),
$$
$$
     |T(t, \omega)(x) | \leq | x | 
e^{- \delta _2 t}\quad  \hbox{for all} \,\,\quad t \geq \tau^* _2, x \in {\Cal S} 
(\omega),
$$
where  $\tau^* _i = \tau^* _i (x,\omega) > 0, i=1,2,$ are  random times 
and $\delta _i > 0, i=1,2,$ are  fixed. 
   
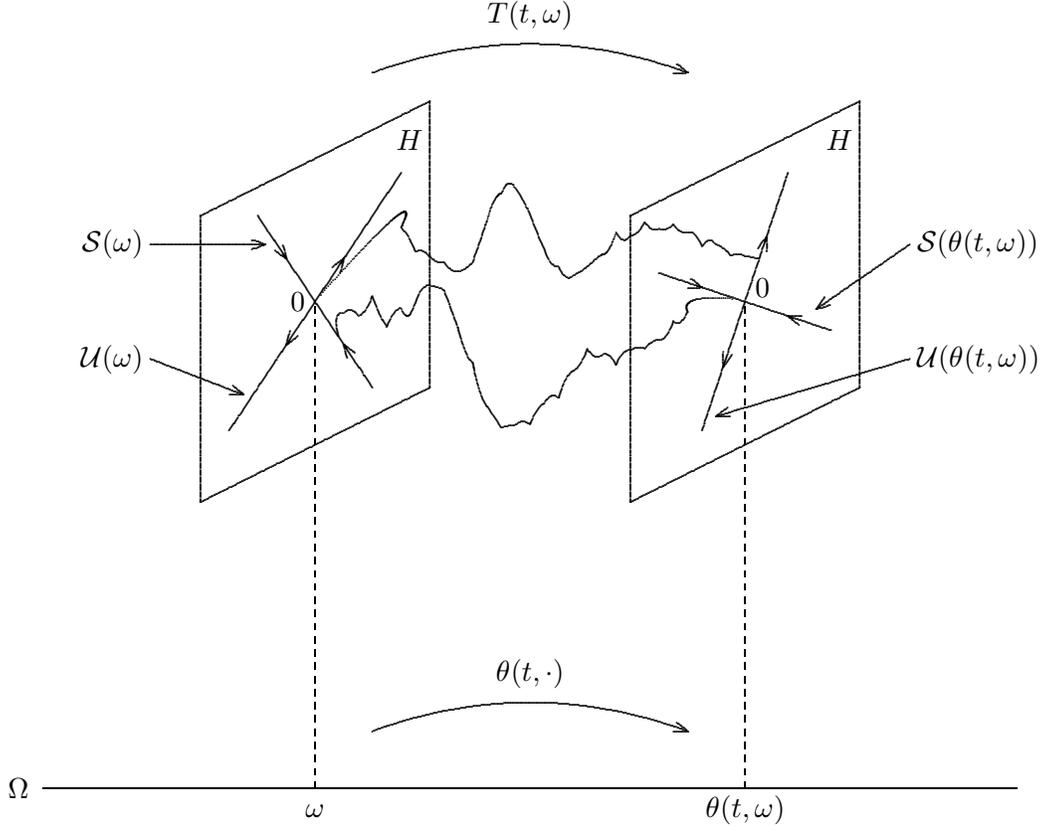
\endproclaim

\bigskip

\centerline{
\beginpicture
\setcoordinatesystem units <0.075truein, 0.075truein>
\setplotarea x from -5 to 65, y from -10 to 50
 %
%
\plot  8 20   24 28   24 48    8 40    8 20 /
\plot 38 20   54 28   54 48   38 40   38 20 /
%
%
\setlinear
\plot 10 25    22 43 /
\plot 12 40    20 28 /
\plot 43 25    49 43 /
\plot 40 36    52 32 /
%
%
\arrow <.07truein> [.4, .7] from    17.8 36.7  to 18 37
\arrow <.07truein> [.4, .7] from    14.2 31.3  to 14 31
\arrow <.07truein> [.4, .7] from    13.8 37.3  to 14 37
\arrow <.07truein> [.4, .7] from    18.2 30.7  to 18 31
\arrow <.07truein> [.4, .7] from 47.4 38.2     to 47.5 38.5
\arrow <.07truein> [.4, .7] from 44.6 29.8     to 44.5 29.5
\arrow <.07truein> [.4, .7] from 42.7 35.1     to 43.0 35.0
\arrow <.07truein> [.4, .7] from 49.3 32.9     to 49.0 33.0
%
%
\putrule from -3 0 to  65 0
\setdashes <2.5truept>
\putrule from 16 0  to  16 33.8
\putrule from 46  0  to  46 34
\put {$\omega$}                 [cb] at 16 -2.0
\put {$\theta(t,\omega)$}     [cb] at 46 -2.5
\put {$\Omega$}                 [rc] at -4 0
\setsolid
\setquadratic    
\plot 20 50   31 52   42 50 /
\arrow <.07truein> [.4, .7] from 41.9 50.034 to 42 50
\put {$T(t,\omega)$} [cb] at 31 53
\plot 20 4    31 6    42 4  /
\arrow <.07truein> [.4, .7] from 41.9 4.034 to 42 4
\put {$\theta(t,\cdot)$} [cb] at 31 7
\put {$\Cal S(\omega)$} [rc] at 4 38
\put {$\Cal U(\omega)$} [rc] at 4 30
\put {$\Cal S(\theta(t,\omega))$} [lc] at 58 38
\put {$\Cal U(\theta(t,\omega))$} [lc] at 58 30
\arrow <.07truein> [.4, .7] from 4.5 38 to 12.5 38
\arrow <.07truein> [.4, .7] from 4.5 30 to 11 27.5
\arrow <.07truein> [.4, .7] from 57.5 38 to 51 33.5
\arrow <.07truein> [.4, .7] from 57.5 30 to 44 26
\put {$H$} [rt] at 23.5 46
\put {$H$} [rt] at 53.5 46
\put {$0$}   [rc] at 15.3 34
\put {$0$}   [lc] at 46.7 35.0
%
%
\setquadratic
\plot
16.4 34.5
21.5 39.8181201
22.0 39.1589922
22.5 38.5453681
23.0 37.4816667
23.5 37.5385127
24.0 37.1480237
24.5 36.8085227
25.0 36.7583332
25.5 36.2309537
26.0 36.1239133
26.5 36.3299163
27.0 36.7666667
27.5 38.2017349
28.0 39.6688230
28.5 41.0514998
29.0 41.5533333
29.5 42.2511695
30.0 41.8882947
30.5 41.0812726
31.0 40.1283333
31.5 38.8092122
32.0 37.6529983
32.5 36.6702851
33.0 36.4633333
33.5 35.7385437
34.0 35.8122121
34.5 36.1047744
35.0 36.4500000
35.5 36.9006751
36.0 37.2856531
36.5 37.6528047
37.0 38.8266667
37.5 38.3228185
38.0 38.6076752
38.5 38.8386944
39.0 39.5666667
39.5 39.0814886
40.0 39.0961458
40.5 39.0627303
41.0 39.3800000
41.5 38.9215400
42.0 38.8202413
42.5 38.6838221
43.0 37.8783333
43.5 38.2636017
44.0 37.9978887
44.5 37.7332315
45.0 37.2200000
45.5 37.3209283
46.0 37.1882038
46.5 37.0863773
47.0 37
/
%
%
\plot
17.6 31.7
17.5 32.8776418
18.0 33.2042269
18.5 33.4286985
19.0 33.0166667
19.5 33.3946939
20.0 34.3998193
20.5 33.1300352
21.0 32.2516667
21.5 33.1920205
22.0 33.5589959
22.5 33.0214732
23.0 33.9766667
23.5 34.9075366
24.0 35.1266973
24.5 35.0325092
25.0 34.7950000
25.5 33.4590829
26.0 32.0592150
26.5 30.5047397
27.0 29.7316667
27.5 27.7170694
28.0 26.6989426
28.5 25.9563444
29.0 25.2200000
29.5 25.3288897
30.0 25.3950149
30.5 25.6386327
31.0 25.4083333
31.5 26.4283096
32.0 26.9084979
32.5 27.4344373
33.0 27.2783332
33.5 28.5906846
34.0 29.1584935
34.5 29.6470559
35.0 30.7200000
35.5 30.1823898
36.0 30.2450278
36.5 30.2601520
37.0 29.5666667
37.5 30.4188190
38.0 30.5988950
38.5 30.8045235
39.0 30.7216667
39.5 31.1657718
40.0 31.3468922
40.5 31.6045666
41.0 32.6166667
41.5 32.5759063
42.0 33.3010359
42.5 34.1256474
45.2 34.2
/
\endpicture
}

\medskip

\medskip


\subheading{2. The non-linear ergodic theorem}

The main objective of this section is to refine and extend discrete-time 
results of D. Ruelle to the continuous-time  setting in Theorem 2.1 below.
This setting  underlies  the dynamics 
of the  semilinear see's and spde's studied by the authors in 
[M-Z-Z.1]. As will be apparent later, the extension of Ruelle's results 
to continuous-time is non-trivial. Indeed, Section 3 in its entirety is devoted 
to the proof of Theorem 2.1. The main difficulties in the analysis are 
outlined after the statement of the theorem.

In the following, denote by $B(x, \rho)$ the open ball, radius $\rho$ and 
 center $x \in H$, and by $\bar B(x, \rho)$ the corresponding closed 
ball.

\bigskip
  
\proclaim{Theorem 2.1} (The local stable manifold theorem)

 Let $(U,\theta)$ be a $C^{k,\e}$ ($k \geq 1, \e \in (0,1]$) perfect cocycle on a 
separable Hilbert space 
$H$ such that for each $(t, \o) \in (0,\infty) \times \O$, $U(t,\cdot,\o): H \to H$ 
takes 
bounded sets into relatively compact sets.
For any $\rho \in (0,\i)$, denote by $\|\cdot\|_{k,\e}$ the $C^{k,\e}$-norm on the 
space 
$C^{k,\e}(\bar B(0,\rho), H)$.  
Let $Y$  be a hyperbolic stationary point of the cocycle 
$(U,\theta)$ satisfying the following integrability property:  
$$ \int_{\O} \log^+ \sup_{0\leq t_1,t_2 \leq a} \|U(t_2, Y(\theta
(t_1,\o)),\theta (t_1,\o))\|_{k,\e} \, dP(\o) < \infty $$
for any fixed $0 <  \rho, a < \infty$  and $\e \in (0,1]$. 
Denote by $\{\cdots < \lambda_{i+1} < \lambda_i < \cdots < \lambda_2 
< \lambda_1 \}$ the Lyapunov 
spectrum of  the linearized cocycle 
$(D U (t, Y(\o),\o), \theta (t,\o), t \geq 0)$.  Define 
$\lambda_{i_0}:=\hbox{max}\{\lambda_i:\lambda_i 
< 0\}$ if at least one $\lambda_i < 0$. If all finite $\lambda_i$ are 
positive, set $\lambda_{i_0} := -\infty$. (Thus 
$\lambda_{i_0-1}$ is the smallest positive Lyapunov exponent of the 
linearized cocycle, if at least one $\lambda_i > 0$; in case all the 
$\lambda_i$'s are negative, set $\lambda_{i_0-1}:=\infty$.)    

Fix $\e_1 \in (0,-\lambda_{i_0})$ and $\e_2 \in  (0,\lambda_{i_0-1})$. 
Then there exist 

\item{(i)} a sure event $\,\O^* \in \F$ with 
$\theta (t,\cdot) (\O^*)= \O^*$ for all $\, t \in \R$, 
\item{(ii)}$\bar \F$-measurable random variables 
$\rho_i,\beta_i: \O^* \to (0, 1), \, \beta_i > \rho_i > 0,\, i=1,2$,    
such that for each $\o \in \O^*$, the following is true:

  There are   
 $C^{k,\e}$ ($\e \in (0,1]$) 
 submanifolds 
 $ \tilde \Cal S (\o),\, \tilde \Cal U (\o)$ of  $\bar B(Y(\o), \rho_1  (\o))$ and 
 \newline $\bar B(Y(\o), \rho_2 (\o))$ (resp.) with the following properties: 
   
\item{(a)} For $\lambda_{i_0} > -\i$, $\tilde \Cal S (\o)$ is the set of all 
$x \in \bar B(Y(\o), \rho_1 (\o))$ such that 
$$|U(n,x, \o)- Y(\theta (n ,\o))| \leq \beta_1 
(\o)\,e^{(\lambda_{i_0}+\e_1)n}$$ 
for all integers $n \geq 0$. If $\lambda_{i_0} =-\i$, then $\tilde \Cal S (\o)$ is the 
set of all 
$x \in \bar B(Y(\o), \rho_1 (\o))$ such that 
$$|U(n,x, \o)- Y(\theta (n,\o))| \leq \beta_1 (\o)\,e^{\lambda n}$$ for all integers 
$n 
\geq 0$ and any $\lambda \in (-\i,0)$.
Furthermore, 
$$\limsup_{t \to \infty} \frac{1}{t} \log |U(t,x, \o)- Y(\theta (t,\o))| 
\leq 
\lambda_{i_0}\tag{2.1}$$ 
for all $x \in \tilde \Cal S (\o)$.
 Each stable subspace $ \Cal S(\o)$  of the linearized cocycle \newline $(D 
U(t,Y(\cdot),\cdot), \theta (t,\cdot))$  is 
 tangent at $Y(\o)$ to the submanifold  
 $\tilde \Cal S (\o)$, viz.  $T_{Y(\o)} \tilde \Cal S (\o)= \Cal S (\o)$. 
In particular, $\hbox{codim}\,\,\tilde \Cal S(\o)= \hbox{codim} \,\, \Cal S (\o),$ is 
fixed and finite.  
     
 \item{(b)} $\displaystyle \limsup_{t \to \infty} \frac{1}{t} \log \biggl [\sup 
\biggl \{ \frac{|U(t,x_1 ,\o)- U(t,x_2,\o)|}{|x_1 
-x_2|}: 
x_1 \neq x_2,\, x_1,$  $x_2 \in \tilde \Cal S (\o) 
\biggr \} \biggr ] \leq \lambda_{i_0}$.  

\newpage

\item{(c)} (Cocycle-invariance of the stable manifolds): 
   
There exists $\tau_1 (\o) \geq 0$ such that 
 $$U(t, \cdot, \o)(\tilde \Cal S  (\o))\subseteq  \tilde 
\Cal S  (\theta (t,\o)) \tag{2.2}$$ 
\indent for all $ t \geq \tau_1 (\o)$. Also
$$
  D U(t,Y(\o), \o) (\Cal S(\o)) \subseteq  \Cal S(\theta (t,\o)),\quad t \geq 0. 
\tag{2.3}
 $$  
 \item{(d)} For $\lambda_{i_0-1} < \i$, $\tilde \Cal U(\o)$ is the set of all
$x \in \bar B(Y(\o), \rho_2 (\o))$ with the property that there is 
a discrete-time
``history" process
$y (\cdot,\o):\{-n: n \geq 0\} \to H$ such that $y(0,\o)=x$ and for each 
integer
$n \geq 1$, one has $U(1, y(-n,\o), \theta (-n,\o))= y(-(n-1),\o)$ and 
$$
|y (-n,\o)-Y(\theta (-n,\o))| \leq \beta_2 (\o) 
e^{-(\lambda_{i_0-1}-\e_2)n}. 
$$ 
 If $\lambda_{i_0-1} =\i$, $\tilde \Cal U(\o)$ is the set of all
$x \in \bar B(Y(\o), \rho_2 (\o))$ with the property that there is 
a discrete-time ``history" process
$y (\cdot,\o):\{-n: n \geq 0\} \to H$ such that $y(0,\o)=x$ and for each
 integer $n \geq 1$,  
$$
|y (-n,\o)-Y(\theta (-n,\o))| \leq \beta_2 (\o) e^{-\lambda n},
$$
for any $\lambda \in (0,\i)$.  
Furthermore, for each $x \in \tilde \Cal U(\o)$, there is 
a unique continuous-time ``history" process also denoted by 
$y (\cdot,\o):(-\infty, 0] \to H$ such that $y(0,\o)=x$, 
$U(t, y(s,\o), \theta (s,\o))= y(t+s,\o)$ for all
$s \leq 0, 0 \leq t \leq -s$, and 
$$
\limsup_{t \to \infty} \frac{1}{t} \log |y(-t, \o)- Y(\theta (-t,\o))| 
\leq -\lambda_{i_0-1}.
$$
Each unstable subspace $\Cal U(\o)$ of the linearized cocycle $(D U(t,Y(\cdot),\cdot), 
\theta (t,\cdot))$   is  
tangent at $Y(\o)$  to $\tilde \Cal U (\o)$, viz. 
$T_{Y(\o)} \tilde \Cal U (\o)= \Cal U (\o)$.
In particular,
$\,\hbox{dim}\,\, \tilde \Cal U (\o)$ is finite and non-random.
 
\item{(e)} Let $y (\cdot,x_i,\o), i =1,2$, be the history 
processes associated with \newline $x_i=y (0,x_i,\o)
 \in \tilde \Cal U(\o), \, i=1,2$.
Then 
$$
\align
\displaystyle \limsup_{t \to \infty} \frac{1}{t} \log \biggl 
[\sup \biggl \{ \frac{|y(-t,x_1 ,\o)- y(-t,x_2,\o)|}
{|x_1 -x_2|}:& 
x_1 \neq x_2,\,
 x_i \in \tilde \Cal U (\o), i=1,2
 \biggr \} \biggr ]\\
& \leq -\lambda_{i_0-1}.
\endalign
$$  
%
 \bigskip
\item{(f)}(Cocycle-invariance of the unstable manifolds): 
   
There exists $\tau_2 (\o) \geq 0$ such that 
 $$\tilde \Cal U  (\o) \subseteq U (t, \cdot, \theta (-t,\o))(\tilde \Cal U (\theta 
(-t,\o)))  
\tag{2.4}
$$ 
\indent
for all  $t \geq \tau_2 (\o)$.
Also
$$
D U(t, \cdot, \theta (-t,\o))( \Cal U (\theta (-t,\o)))=\Cal U (\o), \quad t \geq 
0;
$$ 
\indent and the restriction 
$$
D U (t, \cdot, \theta (-t,\o))| \Cal U (\theta (-t,\o)):  \Cal U (\theta (-t,\o)) 
\to 
\Cal U (\o), \quad t \geq 0,
$$
\indent 
is a linear homeomorphism onto.

\item{(g)} The submanifolds $\tilde \Cal U(\o)$ and $\tilde \Cal S(\o)$ are 
transversal, 
viz.
$$H= T_{Y(\o)} \tilde \Cal U(\o) \oplus T_{Y(\o)} \tilde \Cal S(\o).$$

Assume, in addition, that  the cocycle $(U,\theta)$ is $C^{\i}$. 
Then  the local stable and unstable manifolds  
$\tilde \Cal S(\o),\, \tilde \Cal U(\o)$ are also $C^{\infty}$. 
\endproclaim

\bigskip

The figure below 
summarizes the essential features of the stable 
manifold theorem:

\bigskip

\centerline{\it The Stable Manifold Theorem}

\bigskip
\bigskip

\centerline{
\beginpicture
 \setcoordinatesystem units <0.075truein, 0.075truein>
\setplotarea x from -5 to 65, y from -10 to 50
 %
%
\plot  8 20   24 28   24 48    8 40    8 20 /
\plot 38 20   54 28   54 48   38 40   38 20 /
%
%
\setlinear
\plot 10 25    22 43 /
\plot 12 40    20 28 /
\plot 43 25    49 43 /
\plot 40 36    52 32 /
%
%
\arrow <.07truein> [.4, .7] from    17.8 36.7  to 18 37
\arrow <.07truein> [.4, .7] from    14.2 31.3  to 14 31
\arrow <.07truein> [.4, .7] from    13.8 37.3  to 14 37
\arrow <.07truein> [.4, .7] from    18.2 30.7  to 18 31
\arrow <.07truein> [.4, .7] from 47.4 38.2     to 47.5 38.5
\arrow <.07truein> [.4, .7] from 44.6 29.8     to 44.5 29.5
\arrow <.07truein> [.4, .7] from 42.7 35.1     to 43.0 35.0
\arrow <.07truein> [.4, .7] from 49.3 32.9     to 49.0 33.0
%
%
\putrule from -3 0 to  65 0
\setdashes <2.5truept>
\putrule from 16 0  to  16 33.8
\putrule from 46  0  to  46 34
\put {$\omega$}                 [cb] at 16 -2.0
\put {$\theta(t,\omega)$}     [cb] at 46 -2.5
\put {$\Omega$}                 [rc] at -4 0
\setsolid
\setquadratic    
\plot 16 34     15.2 37  16.6 39.4 /
 \plot 16 34     17 30.5   15.2 29 /
 \plot 20 50   31 52   42 50 /
\plot 50.8 33.9   48 33.7  46 34 /
\plot 46 34  43.8 34.2  40.7 33.0 /
\arrow <.07truein> [.4, .7] from 41.9 50.034 to 42 50
\put {$U(t,\cdot,\omega)$} [cb] at 31 53
\plot 20 4    31 6    42 4  /
\arrow <.07truein> [.4, .7] from 41.9 4.034 to 42 4
\put {$\theta(t,\cdot)$} [cb] at 31 7
%
\put {$\tilde \Cal S(\omega)$} [rc] at 4 36
\put {$\tilde\Cal U(\omega)$} [rc] at 4 30
\put {$\tilde\Cal S(\theta(t,\omega))$} [lc] at 58 38.3
\put {$\tilde\Cal U(\theta(t,\omega))$} [lc] at 58 30
\arrow <.07truein> [.4, .7] from 4.5 36 to 15.5 38
\arrow <.07truein> [.4, .7] from 4.5 30 to 11.5 31.1
 \arrow <.07truein> [.4, .7] from 57.5 37.5 to 50.8 33.9
\arrow <.07truein> [.4, .7] from 57.5 30 to 46.8 31
\put {$ \Cal S(\omega)$} [rc] at 4 40
\arrow <.07truein> [.4, .7] from 4.5 40 to 11.6 40
\put {$\Cal U(\omega)$} [rc] at 4 24.8
\arrow <.07truein> [.4, .7] from 4.5 25 to 10.5 26.2
 \put {$\Cal U(\theta(t,\omega))$} [lc] at 58 46
\arrow <.07truein> [.4, .7] from 57.5 46 to 48.9 41.9
\put {$\Cal S(\theta(t,\omega))$} [lc] at 58 33.9
\arrow <.07truein> [.4, .7] from 57.5 33.5 to 52 31.9
\put {$H$} [rt] at 23.5 46
\put {$H$} [rt] at 53.5 46
\put {$Y(\o)$}   [rc] at 14.5 34
\put {$Y(\theta (t,\o))$}   [lc] at 47.2 35.0
%
\circulararc 360 degrees from 19.47 38 center at 15.76 34.2
\circulararc 360 degrees from 49 38 center at 45.67 34.2
%
\setquadratic
\setdashes <1.5truept>
\plot
20.2 35.6
21.5 37.8181201
22.0 39.1589922
22.5 38.5453681
23.0 37.4816667
23.5 37.5385127
24.0 37.1480237
24.5 36.8085227
25.0 36.7583332
25.5 36.2309537
26.0 36.1239133
26.5 36.3299163
27.0 36.7666667
27.5 38.2017349
28.0 39.6688230
28.5 41.0514998
29.0 41.5533333
29.5 42.2511695
30.0 41.8882947
30.5 41.0812726
31.0 40.1283333
31.5 38.8092122
32.0 37.6529983
32.5 36.6702851
33.0 36.4633333
33.5 35.7385437
34.0 35.8122121
34.5 36.1047744
35.0 36.4500000
35.5 36.9006751
36.0 37.2856531
36.5 37.6528047
37.0 38.8266667
37.5 38.3228185
38.0 38.6076752
38.5 38.8386944
39.0 39.5666667
39.5 39.0814886
40.0 39.0961458
40.5 39.0627303
41.0 39.3800000
41.5 38.9215400
42.0 38.8202413
42.5 38.6838221
43.0 37.8783333
43.5 38.2636017
44.0 37.9978887
44.5 37.7332315
45.0 37.2200000
45.5 37.3209283
45.6 37.1882038
45.6 37.0863773
45.7 36.9
 /
%
%
\setsolid
\plot
17 30.5
17.5 31.0776418
18.0 31.2042269
18.5 32.4286985
19.0 32.9166667
19.5 33.3946939
20.0 33.1998193
20.5 33.0300352
21.0 32.5516667
21.5 33.1920205
22.0 33.5589959
22.5 34.0214732
23.0 34.9766667
23.5 34.9075366
24.0 35.1266973
24.5 35.0325092
25.0 34.7950000
25.5 33.4590829
26.0 32.0592150
26.5 30.5047397
27.0 29.7316667
27.5 27.7170694
28.0 26.6989426
28.5 25.9563444
29.0 25.2200000
29.5 25.3288897
30.0 25.3950149
30.5 25.6386327
31.0 25.4083333
31.5 26.4283096
32.0 26.9084979
32.5 27.4344373
33.0 27.2783332
33.5 28.5906846
34.0 29.1584935
34.5 29.6470559
35.0 30.7200000
35.5 30.1823898
36.0 30.2450278
36.5 30.2601520
37.0 29.5666667
37.5 30.4188190
38.0 30.5988950
38.5 30.8045235
39.0 30.7216667
39.5 31.1657718
40.0 31.3468922
40.5 31.6045666
41.0 32.6166667
41.5 32.5759063
42.0 33.3010359
42.5 33.1256474
43.3 34.0
/
 \setdashes<1.5truept>
\plot 16 34  18.2 35.5  21.0 35 /
\plot  16 34   13 31  11.2 31.5 /
\plot 48 29.8    46.7 31   46 34 /
\plot 46 34  45.7 36.9  44.0 39.0 /
\put{$ t > \tau_1(\o)$} at 31 -7
\endpicture
}

\medskip


\medskip


\baselineskip=24truept

Before we give a detailed proof of Theorem 2.1, we will outline below its 
basic ingredients.

\medskip

\noindent
{\it An outline of the proof of Theorem 2.1:}

\medskip
 
\item{$\bullet$} Since $Y$ is a hyperbolic stationary point of the cocycle 
$(U,\theta)$ (Definition 1.2), then the linearized cocycle satisfies  the hypotheses 
of 
``perfect versions" of 
the ergodic theorem and Kingman's subadditive 
ergodic theorem (Lemma 3.1 (ii), (iii) in Section 3). These refined versions of the 
ergodic theorems give invariance of the
Oseledec spaces under the {\it continuous-time} linearized cocycle (Theorem 1.2). Thus 
the stable/unstable subspaces
will serve as tangent spaces to the local stable/unstable manifolds of the non-linear 
cocycle 
$(U,\theta)$.

\item{$\bullet$} Define the auxiliary perfect cocycle $(Z,\theta)$ by
$$ 
Z(t, \cdot,\o):=U (t, (\cdot)+ Y(\o),\o)- Y(\theta (t,\o)), 
\,\, t \in \R^+, \o \in \O.   
$$
This gives a ``centering" of the cocycle around the stationary trajectory $Y(\theta 
(t))$, with 
the property that $Z$ has a fixed point at $0 \in H$.
Employing the continuous-time integrability estimate in Theorem 2.1, the perfect 
ergodic 
theorem and the perfect subadditive
ergodic theorem, the analysis in ([Ru.2], Theorems 5.1 and 6.1) may be extended to 
obtain local 
stable/unstable manifolds for the
discrete cocycle 
$(Z(n,\cdot,\o),\theta (n,\o))$ near $0$. These manifolds are random objects defined  
for all $\o$ which are sampled from a 
{\it $\theta (t,\cdot)$-invariant sure event} in $\O$.  The translates of these 
manifolds  by the stationary point $Y(\o)$ correspond to 
  local stable/unstable manifolds for $U(n, \cdot,\o)$ near $Y(\o)$.  We then 
interpolate between discrete times and 
extend the arguments in [Ru.2] 
further in order to conclude that the above manifolds for the discrete-time cocycle 
$(U(n, \cdot,\o), \theta (n,\o)), n \geq 1$, also serve   
as local stable/unstable manifolds for the 
{\it continuous-time} cocycle $(U, \theta)$ near $Y$.

 \item{$\bullet$} It turns out that the local stable/unstable manifolds are 
asymptotically 
invariant under the continuous-time cocycle $(U, \theta)$. For the stable manifolds, 
the invariance 
 follows 
by arguments based on (a) a refined version of the perfect subadditive ergodic theorem 
(Lemma 3.2, Section 3), and (b) difficult estimates using the integrability property 
of 
Theroem 2.1 and arguments behind the proofs of Ruelle's Theorems 4.1, 5.1 ([Ru.2]).
 To establish asymptotic invariance of the 
local unstable manifolds, we introduce  the concept of a 
{\it stochastic history process} for $U$, which compensates for the lack of 
invertibility of the cocycle.  Perfection arguments similar  to the 
above give the invariance.   This completes the outline of the proof of 
Theorem 2.1.

\bigskip

A full proof of Theorem 2.1 will be given in the next section. The proof is based on  
a 
discrete-time version of the theorem
given in theorems 5.1, 6.1 
[Ru.2]. The extension to continuous-time  is done via perfection techniques and 
interpolation between discrete times.

\medskip
\subheading{3. Proof of the local stable manifold theorem 
}

The main objective of this section is to give a proof of Theorem 2.1.  In particular, 
we show that the local stable/unstable manifolds for the discrete cocycle 
are parametrized by sure events which are invariant under the continuous-time shift 
$\theta (t,\cdot): \O \to \O$. This is achieved 
via a number of computations based on perfection techniques. Excursions of the cocycle 
between 
discrete times are controlled by integrability hypothesis on the cocycle 
$(U,\theta)$ (Theorem 2.1).

``Perfect versions" of the ergodic theorem and Kingman's subadditive ergodic theorem 
will be used  
to construct the 
shift-invariant sure events appearing in the statement of the local 
stable manifold theorem (Theorem 2.1).  These results are given in Lemmas 3.1 
and 3.2 beow. 
 
 The following convention will be frequently used throughout the paper:

\definition{Definition 3.1}

A family of propositions $\{P(\o): \o \in \O\}$ is said to  
{\it hold perfectly in $\o$\/} if there is a sure event $\O^* \in \F$ such that
$\theta (t,\cdot)(\O^*)=
\O^*$ for all $t \in \R$ and $P(\o)$ is true {\it for every\/} $\o \in \O^*$.
\enddefinition
 
\medskip

\proclaim{Lemma 3.1}

\item{(i)} Let $\O_0 \in \bar \F$ be a sure event such that 
$\theta (t,\cdot)(\O_0) \subseteq \O_0$ for all $t \geq 0$. Then there is
 a sure event $\O_0^* \in \F$ such that $\O_0^*
\subseteq \O_0$ and $\theta (t,\cdot)(\O_0^*)=\O_0^*$ for all $t \in \R$.

\item{(ii)}  Let $h:\O \to \R^+$ be {\it any} function such that there 
exists an $\bar \F$-measurable function $g_1 \in L^1(\O,\R^+;P)$ and a sure event 
$\O_1 
\in \bar \F$ such that 
$\displaystyle \sup_{0 \leq u \leq 1} h(\theta (u,\o)) \leq g_1(\o)$
for all $\,\o \in \O_1$. Then 
 $$\displaystyle \lim_{t \to \infty} \frac{1}{t} h(\theta (t,\o)) =0$$
perfectly in $\o$.
\item{(iii)} Suppose $f : \R^+ \times \O \to \R \cup \{-\infty \}$ is a process such 
that for each $t \in \R^+$, $f(t,\cdot)$ is 
$(\bar \F, \Cal B(\R \cup \{-\infty\}))$-measurable and  the following conditions 
hold:

\itemitem{(a)} 
There is an $\bar \F$-measurable function 
$ g_2\in L^1(\O,\R^+;P)$ and a sure event $\tilde \O_1 \in \bar \F$ such that   
 $\displaystyle \biggl [\sup_{0\leq u \leq 1} f^+(u,\o)+ 
\sup_{0\leq u \leq 1} f^+(1-u,\theta (u,\o))\biggr ] \leq g_2 (\o)$ for all $\,\o \in 
\tilde 
\O_1$. 
  \itemitem{(b)}$f(t_1 +t_2,\o) \leq f(t_1,\o) + f(t_2, \theta (t_1,\o))$ 
for all $t_1, t_2 \geq 0$ and {\bf all} $\o \in \O$. 

\item{} Then there is  
a fixed (non-random) number $f^* \in \R\cup \{-\infty \}$ such that 
$$ \lim_{t \to \infty} \frac{1}{t} f(t,\o) =f^*$$
%
perfectly in $\o$.
\endproclaim

\demo{Proof}

 Assertion (i) is established in Proposition 2.3 
([M-S.3]). 

To prove assertions (ii) and (iii) of the lemma, the reader may adapt the proofs of 
Lemmas 5 and 7 in [Mo.1] and employ assertion (i) above. 
 Cf. also Lemma 3.3 in [M-S.3].
\qed
\enddemo

 \medskip

Lemma 3.2 below is used to construct the {\it continuous-time}    
shift-invariant sure events which appear in the statement of Theorem 2.1. 
In essence, the lemma  is  a continuous-time ``perfect version" of Ruelle's Corollary 
A.2 ([Ru.2], p. 288).

\proclaim{Lemma 3.2}

 Assume that the process $f : \R^+ \times \O \to \R \cup \{-\infty \}$ is 
$(\Cal B(\R^+) \otimes \F, \Cal B(\R \cup \{-\infty\}))$-measurable 
and satisfies the following integrability and subadditivity conditions:

\item{(a)} 
  $\displaystyle \int_{\O} \biggl [\sup_{0\leq t_1,t_2 \leq a} f^+(t_1,\theta 
(t_2,\o))\biggr ]\, dP(\o) < \infty$ for all $\,\, a \in (0,\i)$.
\item{(b)}$f(t_1 +t_2,\o) \leq f(t_1,\o) + f(t_2, \theta (t_1,\o))$ 
for all $t_1, t_2 \geq 0$ and {\bf all} $\o \in \O$.

Then there exists a fixed (non-random) $f^* \in \R \cup \{-\infty \}$ 
 such that the following assertions hold perfectly in $\o$:
\item{(i)}$\displaystyle \lim_{t \to \infty} \frac{1}{t} f(t,\o) =f^*$.   
%
\item{(ii)} Assume $g^* \in \R$ is finite and such that 
$f^* \leq g^*$.  Then for each $\e > 0$, there is an $\bar \F$-measurable 
function $K_{\e}: \O  \to [0, \infty)$ with the following properties   
$$
\align
f(t-s,\theta (s,\o)) &\leq (t-s)g^* + \e t + K_{\e} (\o),\,\,\, 0 \leq s \leq t < 
\infty,\\
K_{\e}(\theta (l,\o)) &\leq K_{\e} (\o)+ \e l, \,\,\, l \in [0,\infty).
\endalign
$$
%
\endproclaim


\demo{Proof}

Applying Lemma 3.1 (iii), it is easy to see that there is an $f^* \in \R\cup 
\{-\infty\}$ such that assertion (i) holds for all $\o$ in a  sure event $\O_2 \in 
\F$ with $\theta (t,\cdot) (\O_2)=
\O_2$ for all $t \in \R$.   
The integrability hypotheses (a) and Lemma 3.1 (i) imply that there is a sure event 
$\O_0 \subseteq \O_2$ such 
that 
$\O_0 \in \F$, $\theta (t,\cdot) (\O_0)= \O_0$  for all $t \in \R$, and 
$\displaystyle 
\sup_{0\leq t_1,t_2 \leq a} f^+(t_1,\theta (t_2,\o)) < \infty$ for all $a \geq 0$ and 
all 
$\o \in \O_0$. Let $g^*$ be a finite number in $[f^*, \i) $.     
 Define the  non-negative process $g:\R^+ \times \O \to \R^+$ by
$$ 
g(t,\o):=\cases\max \{f(t,\o)-tg^*,0\}, \quad &t \geq 0,\, \o \in \O_0,\\
0   &t\geq 0,\, \o \notin \O_0.
\endcases
$$
Then  $g$ is  $(\Cal B(\R^+)\otimes \F,\Cal B(\R^+)) 
$-measurable and satisfies conditions (a) and (b). 
 
Now consider the non-negative process $g' : \R^+ \times \O \to \R^+ $ defined by
$$
g'(t,\o):=\sup_{0\leq s \leq t} [g(s,\o)+ g(t-s,\theta (s,\o))], \quad t \geq 0,\, \o 
\in \O.
$$
Observe that the projection of a 
$(\Cal B(\R^+)\otimes \F)$-measurable set is $\bar \F$-measurable ([Co], p. 281). 
Therefore, $g'$ 
satisfies all the hypotheses of Lemma 3.1 (iii). This gives  a non-negative   
${g'}^* $ such that $\displaystyle \lim_{t \to \infty} \frac{1}{t} g'(t,\o) 
={g'}^*$  
for all $ \o$ in  a sure event $\O_3 \in \F$, with   $\theta (t,\cdot) (\O_3)= 
\O_3$  
for all $t \in \R$.    

We will show next the following convergence in probability:  
$$ \lim_{t \to \infty} \frac{1}{t}\sup_{0\leq s \leq t} g(t-s,\theta (s,\cdot)) =0.   
\tag{3.1}$$
To do this, observe that the process 
$h: \R^+ \times \O \to \R $, $ h(t,\o):= g(t,\theta (-t,\o)),\, t \in \R^+, \,\o\in 
\O$,
satisfies the conditions of Lemma 3.1 (iii). Therefore  
$$ \lim_{t \to \infty} \frac{1}{t} h(t,\cdot) =0$$
almost surely and hence in probability. Pick $\delta, t_0 > 0$  
such that $P( \frac{1}{t} h(t,\cdot) \geq \delta ) < \delta$ for all 
$t \geq t_0$. Let $ t \geq t_0$. Then
$$
\align 
\sup_{0\leq s \leq t} \frac{1}{t}g(t-s,\theta (s,\o)) &\leq   
 \sup_{0\leq s \leq t-t_0} \frac{1}{t}g(t-s,\theta (s,\o)) +   \sup_{t-t_0\leq s \leq 
t} 
\frac{1}{t}g(t-s,\theta (s,\o))\\
&\leq  \sup_{0\leq s \leq t-t_0} \frac{1}{t}g(t-s,\theta (-(t-s),\theta (t,\o)))+ 
\sup_{t-t_0\leq s \leq t} \frac{1}{t}g(t-s,\theta (s,\o)).
\endalign
$$
By condition (a), the second term in the right hand side of the last inequality 
converges to zero in probability. The probability that the first term is less than 
or equal to $\delta$ is  at least $1-\delta$. Hence (3.1) holds.

 It follows easily from (3.1) that $ {g'}^*=0$. This implies that assertion (i)  
 holds  for all $\o$ in a   
sure event $\O_4 \in \F$ with   
$\O_4 \subseteq \O_0 \cap \O_3$ and $\theta (t,\cdot) (\O_4)= \O_4$  for
 all $t \in \R$. To complete the proof of 
assertion (ii), let $\e > 0$ and define  the $(\bar \Cal F, \Cal B(\R^+))$-measurable 
function 
$K_{\e}: \O_4 \to [0, \i)$ by 
$$
K_{\e}(\o):= \sup_{0 \leq s \leq t < \infty} [g(t-s,\theta (s,\o)) -\e t]
$$
for all $\o \in \O_4$.   This completes the proof of the lemma.   \qed

\enddemo

Lemma 3.3 below  is essentially  a ``perfect version'' of Proposition 3.2 in 
[Ru.2], p. 257. Our Lemma 3.2 plays a crucial role in the proof of Lemma 3.3. 
In the statement of the lemma, we will use  $\Cal B_s(L(H))$ to denote the Borel 
$\sigma$-algebra on $L(H)$ generated by the strong topology on $L(H)$, viz. the 
smallest topology on $L(H)$ for which all evaluations $L(H) \ni A \mapsto A(z) \in H, 
z \in H$, are continuous.

\proclaim{Lemma 3.3}

Suppose 
$(T^t (\o), \theta (t,\o)),$  $\, t \geq 0,$ is a perfect cocycle of bounded 
linear operators in $H$ satisfying the following hypotheses:
\item{(i)} The process $\R^+ \times \O \ni (t,\o) \mapsto T^t(\o) \in L(H)$ is $(\Cal 
B(\R^+)\otimes\F, \Cal
B_s(L(H)))$-measurable.
\item{(ii)} The map $\R^+ \times \O \ni (t,\o) \mapsto \theta (t,\o) \in \O$ is 
$(\Cal B(\R^+)\otimes\F, \F)$-measurable, and is a group of ergodic $P$-preserving 
transformations on $(\O, \Cal F, P)$.
\item{(iii)} $\displaystyle E \sup_{0 \leq t_1, t_2 \leq a} \log^+ \|T^{t_2}
(\theta (t_1, \cdot))\|_{L(H)} < \infty$ for any finite $a > 0$. 
\item{(iv)} For each $ t > 0$, 
$T^t (\o)$ is compact,  perfectly in $\o$.
\item{(v)} For  any $u \in H$, the map $[0,\infty)  \ni t \mapsto 
T^t(\o)(u) \in H$ is continuous, perfectly in $\o$.

Let $\{\cdots < \lambda_{i+1} < \lambda_i < \cdots < \lambda_2 < \lambda_1 \}$ be the
Lyapunov spectrum of $(T^t (\o), \theta (t,\o))$, with Oseledec spaces
$$
 \cdots E_{i+1}(\o) \subset E_{i}(\o) \subset \cdots \subset  E_{2}(\o) \subset 
E_{1}(\o)=H.
$$
Let $j_0 \geq 1$ be any fixed integer with $\lambda_{j_0} > -\i$. Let the integer 
function 
$r: \{1,2,\cdots ,Q\} \to \{1,2,\cdots, j_0\}$ ``count'' the multiplicities of 
the Lyapunov exponents in the sense that $r(1)=1,\, r(Q)=j_0$,  and  for each 
$1 \leq i \leq  j_0$, 
the number of integers in $r^{-1}(i)$ is the multiplicity of $\lambda_{i}$.  
Set $V_n (\o):= E_{j_0+1}(\theta (n,\o)), n \geq 0$. 

Then the sequence  $T_n (\o):=T^{1} (\theta ((n-1),\o)), n \geq 1,$ satisfies 
Condition (S) of ([Ru.2], pp. 256-257) perfectly in $\o$ with 
$Q=\hbox{codim}\,  E_{j_0+1}(\o)$. In particular, there is an $\F$-measurable 
set of $Q$ orthonormal
vectors $\{\xi_0^{(1)}(\o), \cdots, \xi_0^{(Q)}(\o)\}$ such that 
$ \xi_0^{(k)}(\o) \in [ E_{r(k)}(\o)\backslash E_{r(k)+1}(\o)]$  for $k=1, 
\cdots, Q$, perfectly in $\o$, and satisfying the following properties: 

Set  
$\xi_t^{(k)}(\o):=\displaystyle\frac{T^{t}(\o)(\xi_0^{(k)}(\o))}{|T^{t}(\o)(\xi_0^{(
k)}(\o))|
}$, 
and for any $u \in H$, write  
$$
u= \sum_{k=1}^Q u_t^{(k)}(\o) \xi_t^{(k)}(\o)   + u_t^{(Q+1)}(\o), \quad 
u_t^{(Q+1)}(\o) 
\in V_0(\theta (t,\o)),\,\, \o \in \O.   
$$
Then for any $\e > 0$, there is an $\bar \Cal F$-measurable  
random constant $D_{\e}(\o) > 0$ such that the following inequalities hold perfectly 
in $\o$:
$$
\aligned
 | u_t^{(k)}(\o)| &\leq D_{\e}(\o) e^{\e t} |u| \\
| u_t^{(Q+1)}(\o)| &\leq D_{\e}(\o) e^{\e t} |u|\\
      D_\e (\theta (l,\o))   &\leq D_{\e}(\o) e^{\e l}
\endaligned 
$$
 for all $t \geq 0, \, 1 \leq k \leq Q$ and for all $l \in [0,\i )$. 

 Furthermore, all the random constants in Ruelle's condition (S) ([Ru.2], pp. 256-257) 
may be chosen to be $\bar \F$-measurable in $\o$.
\endproclaim


\demo{Proof}

Our proof runs along similar lines to that of Proposition 3.2 in [Ru.2]: However, one 
has to  maintain the non-trivial requirement that 
all relevant arguments hold perfectly in $\o$. 
 
It is assumed throughout this proof that the reader is familiar with Ruelle's 
conditions (S): (S1)-(S4) 
as spelled out in ([Ru.2], pp. 256-257). 

Observe first   that $T_n (\o)$ satisfies (S1)  perfectly in $\o$. This holds because 
of  (iii), the perfect cocycle property,
Lemma 3.1  and the proof of Theorem 4 ([Mo.1]). 
Note that, by the ordering of the fixed 
Lyapunov spectrum, relation (3.4) of [Ru.2] holds perfectly.  Denote by $\O^*$ the 
$\theta (t,\cdot)$-invariant sure 
event where (S1) holds. Using ergodicity of $\theta$ and the fact that  
$\hbox{codim}\, V_0(\o)=Q$, for all $\o \in \O^*$, it
follows that  $\hbox{codim}\, V_n(\o)=\hbox{codim}\, E_{j_0+1}(\theta (n,\o))= Q$. 
Therefore,  (S2) 
is satisfied for all $\o \in \O^*$.

We next prove that  (S3) holds perfectly. To do this, we will prove the 
 stronger assertion  
that the {\it continuous-time} cocycle $(T^t (\o), \theta (t,\o))$ satisfies (S3) 
perfectly in $\o$.  Set $\hat 
T^t(\o):= T^t (\o)|V_0(\o), \o \in \O^*, t \geq 0$.
 Hence $\hat T^t(\o)(V_0(\o)) \subseteq V_0(\theta (t,\o))$, and the following 
 cocycle identity 
$$
\hat T^{t_1+t_2}(\o)= \hat T^{t_2}(\theta (t_1,\o))\circ \hat T^{t_1}(\o) 
$$
holds for all $\o \in \O^*, t \geq 0$. Denote 
$F_t(\o):= \log \|\hat T^t(\o)\|,\,\, \o \in \O^*, t \geq 0$. Hypothesis (iii) of the 
lemma easily implies that 
  $\displaystyle E\sup_{0 \leq t_1, t_2 \leq a} F^+_{t_2}(\theta (t_1,\cdot)) < \i$ 
 for any finite $a > 0$.  Furthermore, $(F_t (\o), \theta (t,\o))$ is perfectly 
subadditive because of 
the above cocycle identity.   Applying Lemma 3.1,  we obtain a fixed 
number $F^* \in \R \cup \{-\infty\}$ such that 
$$ \lim_{t \to \infty} \frac{1}{t} F_t(\o)= F^*$$
perfectly in $\o$. Suppose $S=j_0$. When $\lambda_{j_0+1} > -\infty$, set  
$\mu^{(S+1)}:=\lambda_{j_0+1}$; and when   
$\lambda_{j_0+1} = -\infty$, take  
$\mu^{(S+1)}$ to be  any fixed number in 
$(-\infty, \lambda_{j_0})$.  Using (3.5), p. 257 of [Ru.2], it follows that 
$F^* \leq \mu^{(S+1)}$.  
Suppose $\e > 0$ and $ \lambda_{j_0+1} > -\infty$.  Then  by Lemma 3.2(ii), there is   
an $\bar \F$-measurable function $K_{\e}: \O \to [0,\i)$ such that  
$$
\log \|\hat T^{t-s}(\theta (s,\o))\| \leq (t-s) \mu^{(S+1)} + \e t + K_{\e} 
(\o),\quad  0 
\leq s \leq t < \infty,  \tag{3.2}
$$
and   
$$ K_{\e}(\theta (l,\o)) \leq K_{\e} (\o)+ \e l, \quad l \in [0,\infty),$$
 perfectly in $\o$.
When $\lambda_{j_0+1} = -\infty$, the inequality (3.2) holds  where 
$\mu^{(S+1)}$ is 
replaced by any (finite) number in $(-\i, \lambda_{j_0})$. Now let $m, n$ be positive 
integers such that $ m < n$. In (3.2), replace $t$ by $n$ and $s$ by $m+1$ to see 
that $T_n (\o),\, n \geq 1,$ satisfies (S3) perfectly in 
$\o$.   

The rest of this proof will now focus on showing that the sequence $T_n (\o), \, n 
\geq 1,$ also 
satisfies Ruelle's condition (S4)  perfectly in $\o$. Indeed, we will establish the 
stronger statement that  the 
 continuous-time cocycle $(T^t (\o), \theta (t,\o))$ satisfies (S4)  
perfectly in $\o$.  Using the orthogonal decomposition $H = V_0(\theta (t,\o)) \oplus 
V_0(\theta (t,\o))^{\perp}$, 
write     
$$T^t (\o)(\xi)=  \check T^t (\o)(\xi)  +  \tilde T^t (\o)(\xi), \quad  \xi \in H, t 
\geq 0, \o \in \O^*. \tag{3.3}$$
That is,  $\tilde T^t (\o)(\xi) \in V_0(\theta (t,\o))$ and  
$\check T^t (\o)(\xi)  \in V_0(\theta (t,\o))^{\perp}$ are the orthogonal 
projections of  $T^t (\o)(\xi)$ on $V_0(\theta (t,\o))$ and $V_0(\theta 
(t,\o))^{\perp}$, 
respectively. Thus (3.3) defines a family of continuous linear operators 
$\check T^t(\o): H \to V_0(\theta (t,\o))^{\perp} \subseteq H, \,\, 
\tilde T^t(\o): H \to V_0(\theta (t,\o)) \subseteq H, \, t \geq 0$.  We now show   
that 
 the family $(\check T^t (\o), \theta (t,\o)), \o \in \O$, satisfies the perfect 
cocycle 
property in $L(H)$. To prove this, we fix any $\o \in \O$. Then by the cocycle 
property of   $(T^t (\o), \theta (t,\o))$ and (3.3), we obtain
$$
\align
  T^{t_1+t_2}(\o)(\xi)&= T^{t_2}(\theta (t_1,\o))[ T^{t_1}(\o)(\xi)]\\
&=\check T^{t_2}(\theta (t_1,\o))[ \check T^{t_1}(\o)(\xi)]+   
\check T^{t_2}(\theta (t_1,\o))[\tilde T^{t_1}(\o)(\xi)]+ 
\tilde T^{t_2}(\theta (t_1,\o)) [\check T^{t_1}(\o)(\xi)] \\
&\qquad +  \tilde T^{t_2}(\theta (t_1,\o))[\tilde T^{t_1}(\o)(\xi)]. \tag{3.4}
 \endalign
$$
 for all  $t_1, t_2 \geq 0, \xi \in H$. Furthermore, $\check T^{t}(\o)(\xi)=0$ for all
 $\xi \in V_0(\o)$, because $V_0(\o)$ is invariant under the cocycle $(T^t (\o), 
\theta (t,\o))$.
  Thus,  
$\check T^{t_2}(\theta (t_1,\o))[\tilde T^{t_1}(\o)(\xi)]=0$ for all $\xi \in H$, and  
(3.4) yields
$$
\align
 T^{t_1+t_2}&(\o)(\xi) \\
=&\check T^{t_2}(\theta (t_1,\o))[ \check T^{t_1}(\o)(\xi)]+ 
\tilde 
T^{t_2}(\theta (t_1,\o)) [\check T^{t_1}(\o)(\xi)]  +  \tilde T^{t_2}(\theta 
(t_1,\o))[\tilde T^{t_1}(\o)(\xi)].  \tag{3.5}
\endalign
$$
Now
$$
T^{t_1+t_2}(\o)(\xi) 
=\check T^{t_1+t_2}(\o)(\xi)+ \tilde T^{t_1+t_2}(\o)(\xi) \tag{3.6}
$$
for all $\xi \in H$. In the right hand side of (3.5),   
 the first term belongs to  $V_0(\theta (t_1+t_2,\o))^{\perp}$, while  
 the second two terms belong to $V_0(\theta (t_1+t_2,\o))$. So by uniqueness 
of the orthogonal decomposition, it follows from (3.6) and (3.5) that 
$$
\check T^{t_1+t_2}(\o)(\xi)=\check T^{t_2}(\theta (t_1,\o))[ \check 
T^{t_1}(\o)(\xi)]\tag{3.7}
$$
for all $\xi \in H$. Hence $(\check T^t(\o),\theta (t,\o))$ is a 
perfect cocycle  in $L(H)$.
  
 We next verify that both cocycles $(T^t (\o), \theta (t,\o))$  and $ (\check T^t 
(\o), \theta (t,\o))$ 
satisfy the conditions of the perfect Oseledec theorem (Theorem 1.1). To see this,
 note that  
$$E \sup_{0 \leq t_1, t_2 \leq a} \log^+ \|\check T^{t_2}(\theta 
(t_1,\cdot))\|_{L(H)} < \infty \tag{3.8}$$ 
for any finite $a > 0$. This follows immediately from the integrability property 
(iii) of the lemma.  Now apply  Theorem 1.1  to
 $(T^t (\o), \theta (t,\o))$  and $ (\check T^t (\o), \theta (t,\o))$. This gives the 
following limits 
$$
\lim_{t \to \infty} \frac{1}{t} \log  | \check T^{t}(\o)(\xi)|= \check 
l_{\xi},\,\qquad
\lim_{t \to \infty} \frac{1}{t} \log  | T^{t}(\o)(\xi)|=  l_{\xi}
$$
perfectly in $\o$ for all $\xi \in H$, with  $ l_{\xi}, \check l_{\xi}$  fixed numbers 
in $\R \cup \{-\infty\}$. We now apply (3.6) 
  in ([Ru.2], p. 259) to obtain 
$$
\check l_{\xi}= \lim_{n \to \infty} \frac{1}{n} \log  | \check 
T^{n}(\o)(\xi)|=\lim_{n 
\to \infty} \frac{1}{n} \log  | T^{n}(\o)(\xi)|=  l_{\xi}
$$
for a.a. $\o$ and all $\xi \in H\backslash V_0(\o)$. Therefore the equality
$$
\lim_{t \to \infty} \frac{1}{t} \log  | \check T^{t}(\o)(\xi)|=  
\lim_{t \to \infty} \frac{1}{t} \log  | T^{t}(\o)(\xi)|
$$ 
holds perfectly in $\o$ for all $\xi \in H\backslash V_0(\o)$. Hence, relation (3.6) 
in 
([Ru.2], p. 259) may be replaced by the continuous-time ``perfect" relation 
$$
\lim_{t \to \infty} \frac{1}{t} \log  \frac{|\check T^{t}(\o)(\xi)|}{| 
T^{t}(\o)(\xi)|}=0 \tag{3.9}
$$ 
%
for all $\xi \in H\backslash V_0(\o)$.

By ([C-V], Theorem III.6, p. 65) and  Gram-Schmidt orthogonalization, we may select 
a set of 
$Q$, 
$\F$-measurable, orthonormal 
vectors $\{\xi_0^{(1)}(\o), \cdots, \xi_0^{(Q)}(\o)\}$ such that $ 
\xi_0^{(k)}(\o) \in [ E_{r(k)}(\o)\backslash E_{r(k)+1}(\o)] \cap  V_{0}(\o)^{\perp}$ 
 for 
$k=1, \cdots, Q$, perfectly in $\o$.  In  the argument in 
[Ru.2], p. 259, replace (3.6) by (3.9) above, $n$ by $t$, $\xi_n^{(k)}$ by 
$\xi_t^{(k)}(\o):=\displaystyle\frac{T^{t}(\o)(\xi_0^{(k)}(\o))}
{|T^{t}(\o)(\xi_0^{(k)}(\o))|}$, $V_n$ by $V_0(\theta (t,\o))$, 
and $\eta_n^{(k)}$ by $\eta_t^{(k)}(\o):=\displaystyle\frac{\check 
T^{t}(\o)(\xi_0^{(k)}(\o))}{|T^{t}(\o)(\xi_0^{(k)}(\o))|}$. Therefore for $u \in 
H$, we write 
$$
u= \sum_{k=1}^Q u_t^{(k)}(\o) \xi_t^{(k)}(\o)   + u_t^{(Q+1)}(\o), \quad 
u_t^{(Q+1)}(\o) 
\in V_0(\theta (t,\o)),  \tag{3.10}
$$
perfectly in $\o$ for all $t \geq 0$. Furthermore, as in [Ru.2], p. 259,  (3.9) 
implies 
that 
$$
\align 
\lim_{t \to \infty} \frac{1}{t} \log &|\det(\eta_t^{(1)}(\o),\cdots, 
\eta_t^{(Q)}(\o))| 
=0, \tag{3.11}
\endalign
$$
perfectly in $\o$.

It remains to prove that  for each $\e > 0$, there is an $\bar \F$-measurable 
non-negative function 
$D_{\e}: \O \to (0, \infty)$ such that the following inequalities 
$$
\left. \aligned
 | u_t^{(k)}(\o)| &\leq D_{\e}(\o) e^{\e t} |u| \\
| u_t^{(Q+1)}(\o)| &\leq D_{\e}(\o) e^{\e t} |u|\\
      D_\e (\theta (l,\o))   &\leq D_{\e}(\o) e^{\e l}
\endaligned \qquad \right \} \tag{3.12}
$$
hold perfectly in $\o$,  for all $t \geq 0, \, 1 \leq k \leq Q$ and for all $l \in 
[0,\i )$. 

In order to establish the inequalities (3.12), we define 
$$  
D_{\e}(\o):=1 +Q \cdot \sup_{0 \leq s \leq t < \i} e^{-\e 
t}|\det(\eta_{t-s}^{(1)}(\theta(s,\o)),\eta_{t-s}^{(2)}(\theta(s,\o)),\cdots, 
\eta_{t-s}^{(Q)}(\theta (s,\o)))|^{-1} \tag{3.13}
$$
perfectly in $\o$. 

First of all we must show that $D_{\e}(\o)$ is finite  perfectly in $\o$.  
Surprisingly, this 
will require some work. Let $0 \leq s \leq t$. Observe that  
the 
determinant of the linear operator $\check T^{t-s}(\theta (s,\o))$ is given by 
$\displaystyle \frac{|\wedge_{k=1}^Q \check T^{t-s}(\theta 
(s,\o))(v_k)|}{|\wedge_{k=1}^Q v_k|}$ 
for  {\it any choice of basis}  $\{v_1,\cdots,v_Q\}$  in $V_0(\theta (s,\o))^{\perp}$. 
Therefore, the following 
inequalities hold perfectly in $\o$:
$$
\align
|\det(\eta_{t-s}^{(1)}&(\theta(s,\o)),\cdots, \eta_{t-s}^{(Q)}(\theta 
(s,\o)))|^{-1}\\
&=\displaystyle \frac{\Pi_{k=1}^Q |T^{t-s}(\theta (s,\o))(\xi^{(k)}_0(\theta 
(s,\o)))|}{|\det (\check T^{t-s}(\theta (s,\o))(\xi^{(1)}_0(\theta (s,\o))), \cdots, 
\check 
T^{t-s}(\theta (s,\o))(\xi^{(Q)}_0(\theta (s,\o))))|}\\
&=\displaystyle \frac{\Pi_{k=1}^Q [| T^{t-s}(\theta (s,\o))(\xi^{(k)}_0(\theta 
(s,\o)))|] 
\cdot 
|\wedge_{k=1}^Q [\check T^s (\o)(\xi^{(k)}_0(\o))]|}{|\det (\check T^{t-s}(\theta 
(s,\o))(\check T^s (\o)(\xi^{(1)}_0(\o))), \cdots, \check T^{t-s}(\theta 
(s,\o))(\check T^s 
(\o)(\xi^{(Q)}_0(\o))))|}\\
&\leq \displaystyle \frac{\Pi_{k=1}^Q [|T^{t-s}(\theta (s,\o))(\xi^{(k)}_0(\theta 
(s,\o)))|\cdot 
|\check T^s (\o)(\xi^{(k)}_0(\o))|]}{|\det (\check 
T^{t}(\o)(\xi^{(1)}_0(\o)), \cdots, \check T^{t}(\o)(\xi^{(Q)}_0(\o)))|}\\
&=\displaystyle \frac{\Pi_{k=1}^Q [|T^{t-s}(\theta (s,\o))(\xi^{(k)}_0(\theta 
(s,\o)))|\cdot |\check T^s (\o)(\xi^{(k)}_0(\o))|]}{\|
 [\check T^{t}(\o)| V_0(\o)^{\perp}]^{\wedge Q}\|} \tag{3.14}\\
&\leq \displaystyle \frac{\| T^{t-s}(\theta (s,\o))\|^Q \cdot \|\check T^s 
(\o)\|^Q}{\| [\check 
T^{t}(\o)| V_0(\o)^{\perp}]^{\wedge Q}\|}. \tag{3.15}
\endalign
$$
By the integrability condition (iii), it follows  that 
$$
\sup_{0 \leq s \leq t \leq a} \| T^{t-s}(\theta (s,\o))\|^Q \cdot \|\check T^s 
(\o)\|^Q  < \i
$$
perfectly in $\o$ for any finite $a > 0$. 

We now prove that for each finite $a > 0$, 
$$ 
\sup_{0 \leq s \leq t \leq a} |\det(\eta_{t-s}^{(1)}(\theta(s,\o)),\cdots, 
\eta_{t-s}^{(Q)}(\theta (s,\o)))|^{-1} < \i \tag{3.16}
$$
perfectly in $\o$. To see this, define the compact set 
$$
S(\o):=\{(t,v_1,\cdots,v_Q): t \in [0,a], v_k \in V_0(\o)^{\perp}, 
|v_k|=1, <v_k,v_l> =0, 1 \leq k < l \leq Q \}.
$$ 
for $\o \in \O$.  Thus (3.16) will hold if we prove that 
$$ 
\inf_{(t,v_1,\cdots,v_Q) \in S(\o)}|\wedge_{k=1}^Q [\check T^{t}(\o)(v_k)]| > 0 
\tag{3.17}
$$
perfectly in $\o$. 
 
To prove (3.17), we observe that each map $\check T^t(\o)|V_0(\o)^{\perp}: 
V_0(\o)^{\perp} \to V_0(\theta (t,\o))^{\perp}$ is injective for each 
$t \geq 0$ perfectly in $\o$. This is an  easy consequence of the cocycle property and 
the 
fact that $\lambda_{j_0} > -\i$. In fact,
$$|\wedge_{k=1}^Q [\check T^{t}(\o)(v_k)]| > 0  \tag{3.18}$$
for all $(t,v_1,\cdots,v_Q) \in S(\o)$.   Furthermore, the map
$$
[0,a] \times [V_0(\o)^{\perp}]^Q \ni (t,v_1,\cdots,v_Q) \mapsto 
|\wedge_{k=1}^Q [\check T^{t}(\o)(v_k)]| \in [0,\i)
$$
is jointly continuous, by  hypothesis (v) of the lemma. By compactness of $S(\o)$, 
(3.18) implies
 (3.17).  Therefore, (3.16) follows from  (3.15) and (3.17). 

The following convergence
$$ 
\lim_{t \to \i} \frac{1}{t} \log \sup_{0 \leq s \leq t } 
|\det(\eta_{t-s}^{(1)}(\theta(s,\o)),\cdots, \eta_{t-s}^{(Q)}(\theta (s,\o)))|^{-1} 
=0 
\tag{3.19}
$$
holds perfectly in $\o$. To prove this  convergence,  note that (3.14) implies the 
following estimate
$$
\align
|\det(\eta_{t-s}^{(1)}(\theta(s,\o)),\cdots, &\eta_{t-s}^{(Q)}(\theta 
(s,\o)))|^{-1}\\
&\leq \displaystyle \frac{\Pi_{k=1}^Q \{\|[ T^{t-s}(\theta (s,\o))|E_{r(k)}(\theta 
(s,\o)]\| 
\cdot 
\|[\check T^s (\o)| E_{r(k)}(\o)]\|\}}{\| [\check T^{t}(\o)| V_0(\o)^{\perp}]^{\wedge 
Q}\|}
\endalign
$$
for $0 \leq s \leq t$ perfectly in $\o$. Let $\e > 0$ be arbitrary.  Taking 
$\displaystyle \frac{1}{t} \log \sup_{0 \leq s \leq t }$ on both sides of the 
above inequality and applying Lemma 3.2(ii) yields  the following 
inequalities
$$
\allowdisplaybreaks
\align
\frac{1}{t} \log \sup_{0 \leq s \leq t }&|\det(\eta_{t-s}^{(1)}(\theta(s,\o)),\cdots, 
\eta_{t-s}^{(Q)}(\theta (s,\o)))|^{-1}\\
&\leq \frac{1}{t} \sup_{0 \leq s \leq t }\biggl \{ \sum_{k=1}^Q (\log \|[ 
T^{t-s}(\theta 
(s,\o))|E_{r(k)}(\theta (s,\o)]\|+ \log \|[\check T^s (\o)| 
E_{r(k)}(\o)]\|)\biggr \} \\
& \qquad \qquad \qquad \qquad -  \frac{1}{t} \log \|[\check T^{t}(\o)| 
V_0(\o)^{\perp}]^{\wedge Q}\| \\
& \leq \frac{1}{t} \sup_{0 \leq s \leq t }\biggl \{ \sum_{k=1}^Q (t-s)\lambda_{r(k)} 
+ \e t+ K^1_\e(\o)+ \sum_{k=1}^Q s\lambda_{r(k)} + \e s+ K^2_\e (\o) \biggr \} \\
& \qquad \qquad \qquad \qquad -  \frac{1}{t} \log \| [\check T^{t}(\o)| 
V_0(\o)^{\perp}]^{\wedge Q}\|  \\
&= \sum_{k=1}^Q \lambda_{r(k)} + 2 \e +\frac{1}{t}[K^1_\e(\o)+K^2_\e (\o)]
-  \frac{1}{t} \log \| [\check T^{t}(\o)| V_0(\o)^{\perp}]^{\wedge Q}\|, \quad t > 0,
\endalign
$$
perfectly in $\o$, with $K^i_\e (\o), i=1,2,$ finite positive random constants 
independent of $t$.  Therefore, the above inequality implies that 
$$
\align
\limsup_{t \to \i} \frac{1}{t} \log \sup_{0 \leq s \leq t }
|\det(\eta_{t-s}^{(1)}(\theta(s,&\o)),\cdots, \eta_{t-s}^{(Q)}(\theta 
(s,\o)))|^{-1}\\
&\leq \sum_{k=1}^Q \lambda_{r(k)} + 2 \e 
- \liminf_{t \to \i} \frac{1}{t} \log \| [\check T^{t}(\o)| 
V_0(\o)^{\perp}]^{\wedge Q}\|\\
&=\sum_{k=1}^Q \lambda_{r(k)} + 2 \e -\sum_{k=1}^Q \lambda_{r(k)}\\
 &= 2 \e.
 \endalign
$$
Since $\e > 0$ is arbitrary, then 
$$
\limsup_{t \to \i} \frac{1}{t} \log \sup_{0 \leq s \leq t 
}|\det(\eta_{t-s}^{(1)}(\theta(s,\o)),\cdots, \eta_{t-s}^{(Q)}(\theta (s,\o)))|^{-1} 
\leq 0 \tag{3.20}
$$
perfectly in $\o$.  The convergence 
(3.11) immediately implies the inequality
$$
\align
\liminf_{t \to \i} \frac{1}{t} \log \sup_{0 \leq s \leq t 
}&|\det(\eta_{t-s}^{(1)}(\theta(s,\o)),\cdots, \eta_{t-s}^{(Q)}(\theta 
(s,\o)))|^{-1}\\
&\geq \liminf_{t \to \i} \frac{1}{t}\log |\det(\eta_t^{(1)}(\o),\cdots, 
\eta_t^{(Q)}(\o))|^{-1} =0  \tag{3.21}
\endalign
$$
 Thus (3.19) follows from (3.20) and (3.21).  

From (3.16), (3.19) and (3.13), we conclude that  $D_\e (\o)$ is finite 
perfectly in $\o$.

From the definition (3.13) of $D(\o)$, one immediately gets  the last inequality  in 
(3.12).  

It remains to show  the first two inequalities in 
(3.12).  In the right hand side of (3.10),  we look at the terms  
$$ 
\check u(\o)= \sum_{k=1}^Q u_t^{(k)}(\o) \eta_t^{(k)}(\o), \quad u \in H, t \geq 0,
$$
 where $\check u(\o),\, \eta_t^{(k)}(\o), 1 \leq k \leq Q,$ are 
viewed as column vectors in 
$\R^Q$ with respect to the basis $\{ \xi^{(k)}_0(\theta (t,\o)): 1 \leq k \leq Q\}$. 
Using Cramer's rule, the above equation may be solved for  each $u^{(k)}_t(\o)$.  In 
view of (3.13),
this yields the following estimates 
$$
\align
|u^{(k)}_t(\o)|
&= \biggl |\frac{\det(\eta_t^{(1)}(\o),\cdots, \eta_t^{(k-1)}(\o),\check u(\o)
, \eta_t^{(k+1)}(\o),\cdots,\eta_t^{(Q)}(\o))}{\det(\eta_t^{(1)}(\o),\cdots, 
\eta_t^{(Q)}(\o))}\biggr |\\ 
&\leq \frac{|\check 
u(\o)|}{|\det(\eta_t^{(1)}(\o),\cdots, \eta_t^{(Q)}(\o))|}\\ 
&\leq \frac{[D_{\e}(\o)-1]}{Q} |u|e^{\e t} \tag{3.22}\\ 
&\leq D_{\e}(\o) |u|e^{\e t},\quad 1 \leq k \leq Q,\,  t \geq 0, \endalign 
$$ 
perfectly in $\o$. 
By virtue of  (3.10), the triangle inequality and (3.22), one gets 
$$
|u^{(Q+1)}_t(\o)| \leq |u| + \sum_{k=1}^Q |u^{(k)}_t(\o)| \leq D_{\e}(\o) 
|u|e^{\e t}, \quad t \geq 0, 
$$
perfectly in $\o$.  Therefore, $T_n (\o)$ satisfies (S4) perfectly in $\o$, 
and the proof of the proposition is now complete.  \qed
\enddemo

 \medskip


In Lemma 3.4 below, an integrability condition allows us to pass from discrete-time 
limits of the 
cocycle to continuous ones. This property is crucial to the proof of Theroem 2.1.
The reason this property holds is because the integrability hypothesis together with 
the perfect ergodic theorem (Lemma
3.1 (ii)) allow for control of the excursions of the contnuous-time cocycle between 
discrete times.


\proclaim{Lemma 3.4}

 Let $Y: \O \to H$ be a stationary point of the cocycle $(U, \theta)$ satisfying the 
integrability condition
$$ 
\int_{\O} \log^+ \sup_{0\leq t_1,t_2 \leq a} \|U(t_2, Y(\theta
(t_1,\o))+(\cdot),\theta (t_1,\o))\|_{k,\e} \, dP(\o) < \infty $$
for any fixed $0 <  \rho, a < \infty$  and $\e \in (0,1]$. 

Define the random field  $Z: \R^+ \times H \times \O \to H$ by 
$$ Z(t, x,\o):=U(t, x+ Y(\o),\o)- Y(\theta (t,\o)) $$
for $t \geq 0, x \in H, \o \in \O$.  Then $(Z,\theta)$ is a $C^{k,\e}$ perfect 
cocycle. 
Furthermore, there is a sure event
$\O_5 \in  \F$ with the following properties:

\item{(i)}$\theta (t,\cdot)(\O_5)=\O_5$ 
for all $t \in \R$, 
\item{(ii)} For every $\o \in \O_5$ and any $x \in H$, the statement 
$$
\limsup_{n \to \infty} \frac{1}{n } \log |Z(n ,x, \o) | < 0 \tag{3.23}
$$
implies  
$$
\limsup_{t \to \infty} \frac{1}{t} \log |Z(t,x, \o) |= \limsup_{n \to \infty}
 \frac{1}{n } \log |Z(n ,x,\o) |.\tag{3.24} 
$$
 \endproclaim


\demo{Proof}

Note that, by definition, $Z$ is a ``centering'' of the 
cocycle $U$ with respect to the stationary trajectory $\{ Y(\theta (t, \cdot)): t \geq 
0\}$ in the sense that 
$Z(t,0,\o)=0$ for all $(t,\o) \in \R^+ \times \O$. Furthermore, 
$(Z,\theta)$ is a $C^{k,\e}$ perfect cocycle. To see this let $t_1, t_2 \geq 0, \o \in 
\O, x \in H$. Then by the perfect cocycle
property for $U$, it follows that 
$$
\align
Z(t_2, Z(t_1, x,\o), \theta (t_1,\o)) &= 
U(t_2, Z(t_1, x,\o)+Y(\theta (t_1,\o)), \theta (t_1,\o))- Y(\theta (t_2, \theta 
(t_1,\o)))\\
&= U(t_2, U(t_1, x+ Y(\o),\o), \theta (t_1,\o))- Y(\theta (t_2 +t_1,\o))\\
&= Z(t_1+t_2,x, \o).  
\endalign 
$$
Using the integrability condition of the lemma, the proofs of assertions (i) and (ii) 
follow  in the same manner as for
the corresponding assertions in Lemma 3.4 ([M-S.3]). \qed
   
 \enddemo


\medskip


\demo{Proof of Theorem 2.1}


The proof of the theorem consists in two major undertakings:

\item{(a)} Using Ruelle's discrete-time analysis [Ru.2] to show that the assertions of 
Theorem 2.1 
hold for the discretized cocycle, perfectly in $\o$.
\item{(b)} Extending the discrete-time results to continuous time via perfection 
techniques that are
 essentially based on the ergodic theorem and Kingman's 
subaddive ergodic theorem.

Recall the auxiliary cocycle $(Z,\theta)$ defined in Lemma 3.4. 
Consider the random family of maps
$F_{\o}: \bar B(0,1) \to H,\, \o \in \O,$ given by $F_{\o} (x):= Z(1,x, \o), \,\, x 
\in H$, and the time-one shift
$\tau:= \theta (1,\cdot): \O \to \O$. Adopting  Ruelle's notation ([Ru.2], p. 272), we 
set $F^n_{\o}:= F_{\tau^{n-1}(\o)}\circ
\cdots \circ F_{\tau (\o)}\circ F_{\o}$.
Therefore,  $F^n_{\o}= Z(n,\cdot, \o)$ for 
each $n \geq 1$, because $(Z,\theta)$ is a cocycle. By Lemma 3.4, each map $F_{\o}$ is 
$C^{k,\e}$ ($\e \in (0,1]$) on $\bar
B(0,1)$ 
and by the definition of $Z$, it follows that 
$(DF_{\o})(0)= D U(1,Y(\o),\o)$. By the integrability hypothesis of the theorem, it 
is clear that $\log^+  \|D U(1,Y(\cdot),\cdot)\|_{L(H)}$
is integrable. Moreover, in view of the integrability hypothesis on $(U, \theta)$, 
it follows that  the linearized 
continuous-time cocycle $(D U(t, Y(\o),\o), \theta (t,\o))$ and the discrete-time 
cocycle $((DF^n_{\o})(0), \theta (n,\o))$ 
share the same Lyapunov spectrum, viz.:
$$
\{ -\infty < \cdots < \lambda_{i+1} < \lambda_i < \cdots < \lambda_2 < \lambda_1 \}.
$$
(Cf. [Mo.1]). Assume that $\lambda_{i_0}$ is finite, that is  $\lambda_{i_0} \in 
(-\i,0)$. Therefore, under hypotheses (I) of
Theorem 5.1 in ([Ru.2], p. 272), there 
is  a sure event 
$\O_1^* \in \F$ such that $\theta (t,\cdot) (\O_1^*)= \O_1^*$ for all $t \in \R$, 
$\bar \Cal F$-measurable positive random
variables 
$\rho_1,  \beta_1 : \O_1^* \to (0, 1)$, and a random family of 
$C^{k,\e}$ ($k \geq 1, \e \in (0,1]$) stable submanifolds $\tilde \Cal S_d (\o)$ of 
$\bar B(0, \rho_1(\o))$ satisfying the
following properties for each $\o \in \O^*_1$:
$$
\tilde \Cal S_d (\o)= \{ x \in \bar B(0, \rho_1(\o)): |Z(n,x, \o) | \leq \beta_1 (\o) 
e^{(\lambda_{i_0}+\e_1)n} \, \hbox{for all
integers}\, \, n \geq 0\}.\,\, \tag{3.25}
$$
When $\lambda_{i_0}=-\i$, the stable manifold is defined by
$$
\tilde \Cal S_d (\o)= \{ x \in \bar B(0, \rho_1(\o)): |Z(n,x, \o) | \leq \beta_1 (\o) 
e^{\lambda n} \, \hbox{for all
integers}\, \, n \geq 0\}, \,\, \tag{3.25'}
$$
where $\lambda \in (-\i,0)$ is arbitrary. The stable subspace 
$ \Cal S(\o)$  
of the linearized cocycle $(D U(t, Y(\o),\o), \theta (t,\o))$ is tangent to the stable 
manifold $\tilde \Cal S_d (\o)$ at $0$; viz. 
$T_{0} \tilde \Cal S_d(\o)= \Cal S (\o)$.
 In particular,  $\hbox{codim}\,\, \tilde \Cal S_d (\o)$ 
is finite and non-random. 
 Again by  Theorem 5.1 of [Ru.2]), we have the following estimate on the 
 Lyapunov exponent of the Lipschitz 
 constant of $(Z(n,\cdot), \theta (n,\cdot))$ over its stable manifold:
$$
\align
 \limsup_{n \to \infty} \frac{1}{n} \log \biggl 
[\sup_ {x_1,x_2 \in \tilde \Cal S_d (\o)\atop{x_1 \neq x_2}}  \frac{|Z(n,x_1 ,\o)- 
Z(n,x_2,\o)|}{|x_1 -x_2|} \biggr ] \leq
\lambda_{i_0}. \tag{3.26}
\endalign
$$

The statements in the above paragraph hold for all $\o$ in the 
$\theta (t,\cdot)$-invariant sure event 
$\O_1^*$. In order to construct such an  event, we will use perfection arguments 
and the proof of Theorem 5.1 ([Ru.2], p. 272). Assume first that $k=1$ and $\e > 0$. 
Using the notation of [Ru.2], denote 
$T^t (\o):= D Z(t,0,\o),\, f(\o):= \theta (1,\o),\, 
T_n (\o):= D Z(1,0,\theta((n-1),\o))$, for all $\o \in \O$, any positive real $t$ and 
any integer $n \geq 1$.  It is possible to
replace (5.3) in [Ru.2], p. 274) by its continuous-time {\it perfect} analogue
$$ 
\lim_{t \to \infty} \frac{1}{t} \log^+  \|Z(1,\cdot,\theta (t,\o)) \|_{1,\e}=0. 
\tag{3.27}
$$
This follows from   the
integrability hypothesis of the theorem and the perfect ergodic theorem 
(Lemma 3.1 (ii)). More specifically, (3.27) holds for all $\o$ in a sure 
event $\O_1^* \in \F$ such that $\theta (t,\cdot) (\O_1^*)= \O_1^*$ for 
all $t \in \R$.   Assume $\lambda_{i_0} > -\infty$. Adopting the terminology 
of Theorem 1.1 ([Ru.2], p. 248), take 
$S=i_0-1$, and 
$\mu^{(S+1)}=\lambda_{i_0}$. In case   
$\lambda_{i_0} = -\infty$, take  $\mu^{(S+1)}$ to be any fixed number in 
$(-\infty, 0)$. The integrability hypothesis on $U$ coupled with Lemma 3.3 (where $ 
j_0=i_0-1$) imply 
 the existence of  a sure event $\O_2^* \in \F$ such 
that $\O_2^* \subseteq \O^*_1$, $\theta (t,\cdot)(\O_2^*)=\O_2^*$ for all 
$t \in \R$, 
and the sequence $T_n(\o), V_n (\o):= E_{i_0}( \tau^n(\o)), \,\, n \geq 1$, 
satisfies Conditions (S) of [Ru.2], p. 256) for every $\o \in \O_2^*$. 
Pick and fix any  $\o \in \O_2^*$. As in the proof of Theorem 5.1
([Ru.2], pp. 274-278),  the ``perturbation theorem" (Theorem 4.1,
[Ru.2], pp. 262-263) holds for the sequence $T_n(\o), n \geq 1$. Thus the assertions 
in the 
previous paragraph are valid for 
$k=1$ and any $\e \in (0,1]$. When $k > 1$ and $\e \in (0,1])$, we first apply the 
previous 
analysis to the perfect cocycle
$$ \biggl (\check Z (t, x,x_1,\o):=(Z(t,x,\o), 
D Z(t,x,\o)x_1), \theta (t,\o) \biggr ), \quad x, x_1 \in H, \, t \geq 0, 
$$
on $H \oplus H$. Secondly, we use the inductive argument of ([Ru.2], pp. 278-279) to 
show 
 that the $\tilde \Cal S_d(\o)$  are 
$C^{k,\e}$ manifolds $(k > 1, \e \in (0,1])$  {\it perfectly 
in\/} $\o$.    

To establish assertion (a) of the theorem, let  $\tilde \Cal S(\o), \o \in \O_1^*,$  
be the set defined 
therein. Then the definition of $Z$ and property (3.25) of  $\tilde \Cal S_d (\o)$ 
imply 
 that
$$ \tilde \Cal S(\o)=\tilde \Cal S_d (\o)+ Y(\o) \tag{3.28}$$
for all $\o \in \O_1^*$. Thus 
$\tilde \Cal S(\o)$ is a $C^{k,\e}$ manifold $(k > 1, \e \in (0,1])$, with tangent 
space 
  $T_{Y(\o)} \tilde \Cal S (\o)= T_{0} \tilde \Cal S_d (\o)=\Cal S (\o)$. 
In particular,  $\hbox{codim}\,\,\tilde \Cal S(\o)= \hbox{codim} \,\, \Cal S (\o), \o 
\in \O_1^*,$ is 
finite and non-random.

To complete the proof of the inequality (2.1) in part (a) of the theorem, use  (3.26) 
to get 
$$
\limsup_{n  \to \infty} \frac{1}{n } \log |Z(n ,x, \o) | \leq \lambda_{i_0}
$$
perfectly in $\o$ for all $x \in \tilde \Cal S_d (\o)$. In view of Lemma 3.4, we may 
extend the above estimate to cover its continuous-time counterpart. Hence we obtain  
a  sure event $\O^*_3 \subseteq \O^*_2$, $\O_3^* \in \F$,
 such that $\theta (t,\cdot) (\O^*_3)= \O^*_3$ for all $t \in \R$,  and   
$$\limsup_{t  \to \infty} \frac{1}{t } \log |Z(t ,x, \o) | \leq \lambda_{i_0} 
\tag{3.29}$$
for all $\o \in \O^*_3$ and all $x \in \tilde \Cal S_d (\o)$.  The above inequality 
together with 
definition of $Z$ imply the estimate (2.1) of the theorem.

Next, we establish  assertion (b) of the theorem.  To do so, let  $\o \in \O^*_1$ and 
$x \in \tilde \Cal S _d(\o)$. Then by (3.26), it follows that 
there is a positive integer $N_0:=N_0 (\o)$, independent of $x \in \tilde 
\Cal S _d(\o)$, such that $Z(n,x,\o) \in \bar B(0,1)$ for all $n \geq N_0$. Now Lemma 
3.1(ii) gives 
  a $\theta (t,\cdot)$-invariant sure event $\O_3$ such that 
$$
\lim_{t \to \infty} \frac{1}{t}\log^+ 
\sup_{0\leq u \leq 1,\atop (v^*,\eta^*) \in \bar B(0,1)} 
\|D Z(u, (v^*,\eta^*),\theta (t ,\o))\|_{L(H)}=0 
\tag{3.29'}
$$
for all $\o \in  \O_3$.  Define the sure event  $\O^*_4:= \O^*_3 \cap \O_3 \in \F$. 
Clearly,  
   $\theta (t,\cdot)(\O^*_4) =\O^*_4$ for all $t \in \R$. By the definition of $Z$ and 
the Mean Value Theorem, 
we obtain the following inequalities
$$ 
\align
 \sup_{n  \leq t \leq n+1  } \frac{1}{t}  \log \biggl [&\sup_{x_1 \neq 
x_2, \atop (v_1, \eta_1), x_2 \in \tilde \Cal S (\o)} \frac{|U 
(t,x_1, \o)-U (t,x_2,\o)|}{|x_1-x_2|}\biggr ] 
\\    
&= \sup_{n  \leq t \leq n+1  }   \frac{1}{t} \log \, \biggl [\sup_{x_1 
\neq x_2,\atop x_1, x_2 \in \tilde \Cal S_d  (\o) } 
\frac{|Z(t, x_1,\o)- Z(t,x_2,\o)|}{|x_1-x_2|}\biggr ] \\
 &\leq  \, \frac{1}{n }\log^+ \sup_{0\leq u \leq 1,\atop (v^*,\eta^*) \in \bar B(0, 
1)} 
\|D Z(u,(v^*,\eta^*),\theta (n ,\o))\|_{L(H)}
\\
& \qquad \qquad + \frac{n}{(n+1)}\, \frac{1}{n } \log \, \biggl [\sup_{x_1 
\neq x_2,\atop x_1, x_2 \in \tilde \Cal S_d  (\o) } 
\frac{|Z(n, x_1,\o)- Z(n,x_2,\o)|}{|x_1-x_2|}\biggr ]     
\endalign 
$$
for all $\o \in \O^*_4, \,$  and all sufficiently large $\, n \geq N_0(\o)$. 
Now take $\,\displaystyle\limsup_{n \to \infty}\,\,$ on both sides of the above 
inequality,  and 
use (3.26), (3.29$'$) in order to complete the proof assertion (b) of the theorem.

The cocycle invariance (2.3) in part (c) of the theorem follows immediately from  the  
Oseledec-Ruelle theorem 
(Theorem 1.1) applied to the perfect linearized cocycle \newline 
$(D U (t, Y(\o),\o), \theta (t,\o))$. Indeed, one gets  a 
sure $\theta (t,\cdot)$-invariant event  $\in \F$ (also denoted by $\O^*_1$), such 
that $D U (t,Y(\o), \o) (\Cal S(\o)) \subseteq \Cal S(\theta (t,\o))$ for all 
$\, t \geq 0$ and all $\o \in \O^*_1$.

The proof of the asymptotic invariance property (2.2) of the non-linear cocycle 
requires some work.
To achieve this, we will  extend the arguments underlying  the proofs of Theorems 5.1 
and 4.1 in [Ru.2], 
pp. 262-279, to a continuous time setting.  The crucial step towards this goal is to 
show that the  two random variables
$\rho_1, \beta_1$ in (3.25) may be redefined on a sure 
event (also denoted by) $\O^*_1$ such  that $\theta (t,\cdot)(\O^*_1)= 
\O^*_1$ for all $t \in \R$, and 
$$ 
\rho_1 (\theta (t,\o)) \geq  \rho_1 (\o) e^{(\lambda_{i_0}+\e_1) t},
\quad \beta_1 (\theta (t,\o)) \geq \beta_1 (\o)e^{(\lambda_{i_0}+\e_1) t}
\tag{3.30}
$$
for every $\o \in \O^*_1$ and all $t \geq 0$.
  For the given choice of $\e_1$, fix 
$0 < \e_3 < -\e(\lambda_{i_0} + \e_1)/4$, where $\e \in (0,1]$ denotes the H{\"o}lder 
exponent of $U$. 
The above inequalities hold in the {\it discrete\/} 
case (when $t=n$, a positive integer) because of Theorem 5.1 (c)
([Ru.2], p. 274).  To prove them for any continuous time $t$, we will modify the 
definitions of 
 $\rho_1, \beta_1$ in the 
 proofs of Theorems 5.1 and 4.1 in [Ru.2]. In the notation of the proof 
of Theorem 5.1 ([Ru.2], p. 274), we replace the random variable $G$ in 
(5.4) ([Ru.2], p. 274) by the larger one
$$
\tilde G(\o):=\sup_{t \geq 0}\|Z(1,\cdot,\theta (t, \o))\|_{1,\e}\, 
e^{(-t\e_3-\lambda 
\e)}.    \tag{3.31}
$$
Clearly, $\tilde G(\o)$ is finite  perfectly in $\o$, because of (3.27) and Lemma 3.2. 
Following ([Ru.2], pp. 266, 274), the random variables $\rho_1, \beta_1$ may 
be chosen according to the relations
$$
\align 
\beta_1 &:= \biggl [\frac{\delta_1 
\wedge \bigl (\frac{1}{\sqrt{2} A} \bigr )}{2\tilde G}\biggr ]^{\frac{1}{\e}} 
\wedge 1 \tag{3.32}\\  
\rho_1 &:= \frac{\beta_1}{B_{\e_3}}\tag{3.33}\\
\endalign
$$
where $ A, \delta_1$ and $B_{\e_3}$  are random positive constants that are 
 defined via continuous-time analogues of the  
  relations (4.26), (4.18)-(4.21), (4.24), (4.25) in [Ru.2], pp. 265-267, 
with  $\eta$  replaced by $\e_3$. In particular, the ``ancestry'' of $ A, \delta_1$ 
and 
$B_{\e_3}$ in Ruelle's argument may be traced back to the constants  
    $D_{\e_3}, K_{\e_3}$ which appear 
in Lemmas 3.3 and 3.2 of this article. Hence (3.30) will follow if we can show that, 
for sufficiently small $\e_3 > 0$, the following inequalities
$$
\left. \aligned
 K_{\e_3}(\theta(l,\o)) &\leq K_{\e_3}(\o) + {\frac{\e_3 l}{2}} \\
  D_{\e_3}(\theta(l,\o)) &\leq e^{\frac{\e_3 l}{2}}D_{\e_3}(\o)\\
  \tilde G(\theta (l,\o)) &\leq e^{\e_3 l} \tilde G(\o)
  \endaligned \right \} \tag{3.34}
 $$
hold perfectly in $\o$ for all real $l \geq 0$. The first and second inequalities in 
(3.34)
 follow from Lemmas 3.2(ii) and 3.3, respectively.
The third inequality is an immediate consequence of the definition of $\tilde G$ in  
(3.31). 
The proof of (3.30) is now complete in view of (3.32), 
(3.33) and (3.34).   

The inequalities in (3.30) will allow us to establish the asymptotic invariance 
property (2.2) in (c) of 
the theorem.  By the perfect inequality in (b), there is a sure event $ \O^*_5 
\subseteq \O^*_4$ 
such that $\theta (t,\cdot)(\O^*_5)= \O^*_5$ for 
all $t \in \R$, and  for any $0 < \e' < \e_1$ and  any $\o \in  \O^*_5$, there exists  
$\beta^{\e'}(\o)> 0$ (independent of $x$) so that   
$$
|U (t,x,\o)-Y(\theta (t,\o))| \leq   \beta^{\e'}(\o)e^{(\lambda_{i_0} 
+\e')t} \tag{3.35}
$$ 
for all $x \in \tilde \Cal S  (\o),  \, t \geq 0$.  Let $t$ be any positive real, 
$n$ a non-negative integer, 
$\o \in  \O^*_5$ and $x \in \tilde \Cal S  (\o)$.  Using the cocycle property and 
(3.35), we 
obtain   
$$
\align 
|U (n,U (t,x,\o),\theta (t,\o))-Y(\theta (n,\theta (t,\o)))|
&=|U (n+t,x,\o) -Y(\theta (n+t,\o))|\\
&\leq \beta^{\e'}(\o)e^{(\lambda_{i_0} +\e')(n+t)} \\
&\leq \beta^{\e'}(\o)e^{(\lambda_{i_0} +\e')t} e^{(\lambda_{i_0}+\e_1)n}.
\tag{3.36}
\endalign
$$
Using  (3.30),(3.35), (3.36) and the 
definition of $\tilde \Cal S  (\theta (t,\o))$, we see that for each  $\o \in \O^*_5$,
 there exists $\tau_1 (\o) > 0$ such that $U(t, x,\o) \in \tilde \Cal S  (\theta 
(t,\o))$ 
for all $t \geq \tau_1 (\o)$. Hence, for all  $\o \in \O^*_5$, 
$$U(t, \cdot, \o)(\tilde \Cal S  (\o))\subseteq  \tilde 
\Cal S  (\theta (t,\o)), \quad t \geq \tau_1 (\o) $$ 
 and the proof of assertion (c) is complete.

Our next objective is to establish the  existence of the {\it perfect} family of
local unstable manifolds $\tilde \Cal U(\o)$ in  assertion (d) of the theorem. 
 To this end, we  define the random field 
$\hat Z: \R^+ \times H \times \O \to H$ by
$$
\hat Z(t, x,\o):= U(t,x+ Y(\theta (-t,\o)), \theta (-t,\o))-Y(\o) 
\tag{3.37}
 $$ 
for all $t \geq 0,\,\,x \in H,\,\, \o \in \O$. Note that 
$\hat Z(t, \cdot,\o)= Z(t,\cdot,\theta(-t,\o)), \, t \geq 0,\, \o \in \O$; 
and $\hat Z$ is 
$(\Cal B(\R^+)\otimes \Cal B(H)\otimes \F, \Cal B(H))$-measurable.
 Since $Y$ is a stationary point for 
$(U, \theta)$, we may replace $\o$  by $\theta (-t,\o)$ in (1.1). Thus  
 $ \hat Z(t, 0,\o) = 0$ for all $t \geq 0,\, \o \in \O$. We contend that 
 $([D \hat Z (t,0,\o)]^*, \theta (-t,\o), \, t \geq 0)$ is a 
perfect linear cocycle in $L(H)$. To see this,  we first observe  that 
$(D U(t, Y(\o),\o), \theta (t,\o))$ is an
$L(H)$-valued  perfect cocycle:
$$ 
D U(t_1+t_2, Y(\o), \o) =D U(t_1, Y(\theta (t_2,\o)), \theta (t_2, \o))\circ D U(t_2, 
Y(\o),\o) 
$$
for all $\o \in \O, t_1, t_2 \geq 0$. Secondly, we replace $\o$ by $\theta 
(-t_1-t_2,\o)$ 
 and take adjoints in the above identity to obtain  
$$ 
\align
[D& U(t_1+t_2, Y(\theta (-t_1-t_2,\o)), \theta (-t_1-t_2,\o))]^* \\
&=[D U(t_2, Y(\theta (-t_1-t_2,\o)), \theta (-t_1-t_2, \o))]^*\circ [D U(t_1, Y(\theta 
(-t_1,\o)),\theta (-t_1,\o)]^* 
\endalign
$$
for all $\o \in \O, t_1, t_2 \geq 0$. Therefore,
$$ 
[D \hat Z(t_1+t_2,0, \o)]^* =[D \hat Z(t_2,0, \theta (-t_1, \o))]^*\circ [D \hat 
Z(t_1,0,\o)]^* 
$$
for all $\o \in \O, t_1, t_2 \geq 0$; and our contention is proved.
 
  We will now apply the Oseledec-Ruelle theorem to the perfect cocycle 
$([D \hat Z (t,0,\o)]^*,\quad$  $\theta (-t,\o), \, t \geq 0)$. To do this, it is 
sufficient to check    the 
integrability condition
$$
 \int_{\O} \log^+ \sup_{0\leq t_1,t_2 \leq a} 
\|[D \hat Z(t_2,0, \theta (-t_1,\o))]^*\|_{L(H)} \, dP(\o) < \infty 
\tag{3.38} $$
for any fixed $a \in (0,\i)$.  The above integrability relation  follows  from the 
integrability hypothesis 
of Theorem 2.1 and the $P$-preserving property of  $\theta (t,\cdot)$:  
$$
\allowdisplaybreaks
 \align
\int_{\O} &\log^+ \sup_{0\leq t_1,t_2 \leq a} \|[D \hat Z(t_2,0, \theta 
(-t_1,\o))]^*\|_{L(H)} \, dP(\o)\\
&=\int_{\O} \log^+ \sup_{0\leq t_1, t_2 \leq a} 
\|D U(t_2, Y(\theta (-t_2-t_1,\o)), 
\theta (-t_2-t_1,\o)))\|_{L(H)} \, dP(\o)\\
&\leq \int_{\O} \log^+ \sup_{0\leq t_1 \leq 2a,\, 0 \leq t_2 \leq a } 
\|D U(t_2, Y(\theta (t_1,\o)), \theta (t_1,\o)))\|_{L(H)} \, dP(\o)\\
&\leq \int_{\O} \log^+ \sup_{0\leq t_1 \leq a,\, 0 \leq t_2 \leq a } 
\|D U(t_2, Y(\theta (t_1,\o)), \theta (t_1,\o)))\|_{L(H)} \, dP(\o)\\
&+\int_{\O} \log^+ \sup_{a\leq t_1 \leq 2a,\, 0 \leq t_2 \leq a } 
\|D U(t_2, Y(\theta (t_1-a,\o)), \theta (t_1-a,\o)))\|_{L(H)} \, dP(\o)\\
&=2\int_{\O} \log^+ \sup_{0\leq t_1, t_2 \leq a } 
\|D U(t_2, Y(\theta (t_1,\o)), \theta (t_1,\o)))\|_{L(H)} \, dP(\o) < \infty.
\endalign
$$
 By (3.38) and the Oseledec-Ruelle theorem, we conclude that the linear 
cocycle  \newline $([D \hat Z (t,0,\o)]^*, \theta (-t,\o), \, t \geq 0)$  has a fixed 
discrete Lyapunov spectrum. Furthermore, this spectrum (with multiplicities) 
coincides with that of  
the cocycle $(D U(t, Y(\o),\o), \theta (t,\o))$, viz. 
$\{\cdots \lambda_{i+1} < \lambda_i < \cdots < \lambda_2 < \lambda_1 \}$ 
where $\lambda_i \neq 0$ for all $i \geq 1$, by hyperbolicity. See  [Ru.2], 
Section 3.5, p. 261. 

The next step in our construction  of the {\it perfect}  random family of 
local unstable manifolds $\tilde \Cal U  (\o)$ is the following estimate:  
$$
\int_{\O} \log^+ \sup_{0\leq t_1,t_2 \leq 1} \|\hat Z(t_2,\cdot, \theta 
(-t_1,\o))\|_{k,\e} \, dP(\o) < \infty.
$$
By the same argument as in the previous paragraph, the above estimate is a consequence 
of  the $P$-preserving 
property of $\theta (t,\cdot), t \in \R$, and the integrability hypothesis 
of the theorem.
Define  $\lambda_{i_0-1}$ as in the statement of Theorem 2.1, and fix any $\e_2 \in 
(0, \lambda_{i_0-1})$.
 In view of the 
above integrability property, it follows from Lemma 3.3 that the sequence 
$\tilde T_n (\o):=[D \hat Z(1,0, \theta (-n, \o))]^*$, $ \, \theta (-n,\o)$, 
$n \geq 0$, 
satisfies Condition (S) of [Ru.2] perfectly in $\o$. Hence the sequence 
$\tilde T_n (\o), n \geq 1,$ satisfies
Corollary 3.4 ([Ru.2], p. 260) perfectly in $\o$, because of Proposition 3.3 
in [Ru.2].  At this point, we may modify the arguments in the proof of Ruelle's 
Theorem 6.1 
([Ru.2], p. 280) using an approach analogous to the one used in constructing the 
stable 
manifolds in this proof. Therefore, one gets  a
$\theta (-t,\cdot)$-invariant sure event $\hat\O_1^* \in \F$ and 
$\bar \Cal F$-measurable random variables $\rho_2, \beta_2:\hat\O_1^* \to (0,1)$
satisfying the following properties.  If $\lambda_{i_0-1} < \i$, define
$\tilde \Cal U_d  (\o)$ to be the set of 
all $x_0 \in \bar B(0, \rho_2 (\o))$ with the property that there is 
a discrete ``history" process $u(-n,\cdot):\O \to H, n \geq 0$, 
such that $u(0,\o)=x_0$, 
$\hat Z (1,u(-(n+1),\o), \theta (-n, \o))=u(-n,\o)$ and 
$|u(-n,\o)| \leq \beta_2(\o) e^{-n(\lambda_{i_0-1}-\e_2)}$ for all 
$n \geq 0$.  If $\lambda_{i_0-1}=\i$, define $\tilde \Cal U_{d}(\o)$ to be the set 
of all $x_0 \in H$ with the property that there is a discrete 
history process $u(-n,\cdot):\O \to H, n \geq 0$, such that 
$u(0,\o)=x_0$, and 
$|u(-n,\o)| \leq \beta_2(\o) e^{-\lambda n}$ for all 
$n \geq 0$ and arbitrary $\lambda > 0$.
It follows from ([Ru.2], p. 281) that the discrete history process $u(-n,\cdot)$ is 
uniquely determined by 
$x_0$.  Moreover, each  
$\,\tilde \Cal U_d  (\o), \o \in \hat\O_1^*$,  is a   $C^{k,\e}$ 
($k \geq 1, \e \in (0,1]$)  finite-dimensional
submanifold of $\bar B(0, \rho_2 (\o))$ with tangent space $\Cal U(\o)$ at $0$, and 
$\hbox{dim}\, \tilde \Cal U_d (\o)$ is fixed independently of 
$\o$ and $\e_2$.  Furthermore, 
$$ \rho_2 (\theta (-t,\o)) \geq  \rho_2 (\o) e^{-(\lambda_{i_0-1}-\e_2) t},
\quad \beta_2 (\theta (-t,\o)) \geq \beta_2 (\o)e^{-(\lambda_{i_0-1}-\e_2) t}.
\tag{3.39}$$
perfectly in $\o$ for all $t \geq 0$.
We claim that the set $\tilde \Cal U  (\o)$ defined in (d) of Theorem 2.1
 coincides with  
$\tilde \Cal U_d  (\o)+ Y(\o)$ for each $\o \in \hat\O^*_1$. 
  We first show that 
$\tilde \Cal U _d (\o)+ Y(\o) \subseteq \tilde \Cal U (\o)$. Let 
$x_0 \in \tilde \Cal U _d 
(\o)$ and $u$ be as above. Set
$$ y_0(-n,\o):= u(-n,\o)+ Y(\theta (-n,\o)), \qquad n \geq 0. \tag{3.40} $$
 It is easy to check that $y_0$ is a discrete history process satisfying the 
first and 
second assertions in (d) of the theorem. Hence 
$x_0+Y(\o) \in \tilde \Cal U (\o)$.    Similarly,   
$\tilde \Cal U (\o) \subseteq  \tilde \Cal U _d (\o)+ Y(\o)$ 
for all $\o \in \hat\O_1^*$. Hence 
$\tilde \Cal U (\o) =  \tilde \Cal U _d (\o)+ Y(\o)$ 
for all $\o \in \hat\O_1^*$. This immediately implies that $\tilde \Cal U  (\o)$ is a   
$C^{k,\e}$ ($k \geq 1, \e \in (0,1]$) 
finite-dimensional submanifold of $\bar B(Y(\o), \rho_2 (\o))$, and 
$$ 
T_{Y(\o)} \tilde \Cal U  (\o) = T_{0} \tilde \Cal U_d  (\o)=  \Cal U(\o)
$$
for all  $\o \in \hat\O_1^*$.

We will next address the issue of the existence of the continuous-time history process 
satisfying  the third assertion in part (d) of the theorem. Suppose 
$x \in \tilde \Cal U (\o)$. From what we proved in the previous paragraph,  it follows 
that 
 there is an $x_0 \in \Cal U _d (\o)$ such that $x=x_0+Y(\o)$. The discrete process 
$y_0$ given by (3.40) may be extended to  
a continuous-time history process $y (\cdot,\o):(-\infty, 
0] \to H$ such that $y(0,\o)=x$, and $y(\cdot,\o)$  satisfies the third  
assertion in (d).  This is achieved by  interpolation within the periods  $[-(n+1), 
-n],\,\, n \geq 0,$ using 
 the cocycle property of  $U$: Indeed, let $s \in  (-(n+1), -n)$. Then there is an  
$\alpha \in (0,1)$, such that
 $s = \alpha -(n+1)$.  Define 
$$ y(s, \o):= U(s+n+1, y_0(-(n+1), \o), \theta (-(n+1), \o)). $$
Obviously,  $y(0,\o)=x_0+Y(\o)=x$. 
 Let $s \in  (-(n+1), -n)$  and suppose $0 < t \leq -s$. Pick  a 
positive integer $m < n$ such that 
$s+t \in [-(m+1), -m]$. The above definition of $y$, together with the perfect cocycle 
property for $U$, 
easily imply that 
$$ 
\align
y(t+s, \o)=U(t, y(s,\o), \theta (s,\o)). \tag{3.41}
\endalign
$$
In particular,  $U(t, y(-t,\o), \theta (-t,\o))= x$ for all $t \geq 0$.  This follows 
from (3.41) when $s$ is replaced 
by $-t$. Furthermore, for each $x \in \tilde \Cal U(\o)$, the above continuous-time 
history process is 
uniquely determined because its discrete-time counterpart is unique.  

   We will now prove the following estimate 
$$
\limsup_{t \to \infty} \frac{1}{t} \log |y(-t, \o)- Y(\theta (-t,\o))| 
\leq -\lambda_{i_0-1}  \tag{3.42}
$$
perfectly in $\o$. We start with its discrete-time counterpart   
$$
\limsup_{n \to \infty} \frac{1}{n} \log |y(-n, \o)-Y(\theta (-n,\o)) |
 \leq -\lambda_{i_0-1}  \tag{3.43}
$$
which holds perfectly in $\o$, because of Theorem 6.1 (b) in [Ru.2].  Let $t \in (n, 
n+1)$. Then  there exists 
$\gamma \in (0,1)$ such that $-t=\gamma -(n+1)$.  Thus by the definition 
of $y$ and the Mean Value Theorem, it follows that
%
$$
\align
|&y(-t, \o)-Y(\theta (-t,\o)) |\\
&=|U(\gamma, y(-(n+1), \o), \theta (-(n+1), \o))- U(\gamma, Y(\theta (-(n+1),\o), 
\theta 
(-(n+1),\o))|\\
&\leq \sup_{ (v^*,\eta^*)\in \bar B(0,1), \atop{\gamma \in (0,1)}}
\|D U(\gamma, (v^*,\eta^*)+Y(\theta(-(n+1),\o)),\theta (-(n+1),\o))\|_{L(H)}\\
&\qquad \quad \times |y(-(n+1), \o)- Y(\theta (-(n+1),\o)))|
\endalign
$$
perfectly in $\o$. Letting $t \to \infty$, we get
$$
\align
&\limsup_{t \to \infty} \frac{1}{t} \log |y(-t, \o)-Y(\theta (-t,\o)) |\\
&\leq \limsup_{n \to \infty} \frac{1}{n} \log^+  \sup_{ (v^*,\eta^*)\in \bar B(0,1), 
\atop{\gamma \in (0,1)}}
\|D U(\gamma, (v^*,\eta^*)+Y(\theta(-(n+1),\o)),\theta (-(n+1),\o))\|_{L(H)}\\
&\qquad + \limsup_{n \to \infty} \frac{1}{n} \log |y(-(n+1), \o)- Y(\theta 
(-(n+1),\o)))|. 
 \endalign
$$
By the integrability condition of the theorem and the perfect ergodic theorem (Lemma 
3.1 (ii)), 
the first term on the right hand side of the above inequality is zero, perfectly in 
$\o \in \O$.  Since $y(0) \in \tilde \Cal U (\o)$, the second term is less than or 
equal 
to $-\lambda_{i_0-1}$.  This completes the proof of assertion (d) of the theorem.

We will omit the proof of  assertion (e), since it is very similar to that of (3.42).

Our next objective is to prove assertion (f) of the theorem. Note first that the 
perfect invariance    
$$
D U(t, \cdot, \theta (-t,\o))( \Cal U (\theta (-t,\o)))=\Cal U (\o), \quad t \geq 
0,
$$ 
%
follows from  the cocycle property for the linearized
semiflow and Theorem (1.2); cf.  [Mo.1], Corollary 2 (v) of Theorem 4. Since 
$\hbox{dim}\,\Cal U (\o)$ is fixed and finite
perfectly in $\o$,   the restriction 
$$
D U (t, \cdot, \theta (-t,\o))| \Cal U (\theta (-t,\o)):  \Cal U (\theta (-t,\o)) 
\to 
\Cal U (\o), \quad t \geq 0,
$$
is a linear homeomorphism onto. It remains to check 
the following asymptotic invariance property in (f): 
$$
\tilde \Cal U  (\o) \subseteq U(t, \cdot, \theta (-t,\o)) (\tilde \Cal U  (\theta 
(-t,\o))), \quad t \geq \tau_2(\o),  
  \tag{3.44}
$$
perfectly in $\o$ for some $\tau_2 (\o) > 0$. Suppose  
 $x \in \tilde \Cal U (\o)$.  Then by assertions (d), (e) of the theorem and 
inequalities (3.39), there 
exist a (unique) history process $y(-t,\o), t \geq 0$, and a random time 
$\tau_2 (\o) > 0$ satisfying the following: $y(0,\o)=x$, 
$y(-t,\o) \in \bar B(Y(\theta (-t,\o)),\rho_2(\theta (-t,\o)))$ 
for all $t \geq \tau_2(\o)$, and 
$$
y(t'-t, \o) =U(t', y(-t,\o), \theta (-t,\o)), \quad 0 <  t' \leq  t, \tag{3.45}
$$
perfectly in $\o$. Pick any  $t_1 \geq \tau_2 (\o)$. Then
  $x=  U(t_1, y(-t_1,\o), \theta (-t_1,\o))$, because of (3.45) (for $t=t'=t_1$). Now  
$y(-t_1,\o) \in \tilde \Cal U(\theta (-t_1,\o)))$. 
To prove this, we define the process 
$y_1(-t,\o):= y(-t-t_1,\o),$  $ t \geq 0$. Hence $y_1 (\cdot,\o)$ is a history process 
and 
$$
y_1 (0,\o)=y(-t_1,\o) \in \bar B(Y(\theta (-t_1,\o)),\rho_2(\theta (-t_1,\o))).
$$
Therefore $y(-t_1,\o) \in \tilde \Cal U(\theta (-t_1,\o)))$. This implies (3.44) 
because  $t_1 \geq \tau_2(\o)$ is arbitrary.

To prove the transversality property in in (g), note the following perfect identities:  
$$
T_{Y(\o)} \tilde \Cal U (\o)= \Cal U (\o), \quad  T_{Y(\o)} \tilde \Cal S (\o)= \Cal 
S 
(\o), \quad  H= \Cal U (\o) \oplus \Cal S (\o).
$$

 All the assertions (a)-(g) of the theorem will hold perfectly in $\o$ if we take   
$\O^*:= \O^*_1 \cap \hat \O^*_1$.
  
To deal with the case when $U$ is a $C^{\i}$ cocycle, we  adapt 
the proof  in [Ru.2], section (5.3) (p. 297). Thus we  obtain a 
$\theta (t,\cdot)$-invariant sure event in 
$\F$ (also denoted by $\O^*$) such that $\tilde \Cal S(\o)$ and $ \tilde \Cal
U(\o)$ are $C^{\infty}$ for all $\o \in \O^*$. This completes  the proof of 
Theorem 2.1. 
\qed    
\enddemo

\subheading{4. The local stable manifold theorem for see's and spde's}

In this section, we discuss several classes of semilinear stochastic evolutions 
equations and spde's. The objective is 
to establish sufficient conditions for a local stable manifold theorem for each class. 

\medskip
\noindent
{\it (a) Stochastic semilinear evolution equations: Additive noise.}

Let $K, H$ be two separable real Hilbert spaces.  Let $A$ be a
self-adjoint operator on $H$ such that $A\geq cI_H$, where $c$ is a real  constant and 
$I_H$ is  the
identity operator on $H$. Assume that $A$ admits
a discrete non-vanishing spectrum $\{ \mu_n, n \geq 1 \}$
which is bounded below. Let $\{e_n,\, n\geq 1\}$ denote a basis for $H$ consisting of 
eigen vectors of $A$,
viz.  $Ae_n=\mu_n e_n, \, n \geq 1$. Assume further that $A^{-1}$ is trace-class.
Suppose $B_0 \in L_2(K,H)$. Let $W(t),  t\in \R,$ be cylindrical  Brownian
motion on the canonical filtered Wiener space $(\O,\F,(\F_t)_{t \geq 0},P)$ and with a
separable covariance Hilbert space $K$ ([M.Z.Z.1], section 2). Let $T_t=e^{-At}$ stand 
for the strongly 
continuous semigroup generated by $-A$. 
%
%
%
%
 
 Denote by  $\mu_m$ the largest negative eigenvalue of $A$ and
by $\mu_{m+1}$ its smallest positive eigenvalue. Thus there is an  orthogonal 
$\{T_t\}_{t \geq 0}$-invariant splitting of $H$ using the negative  eigenvalues 
$\{\mu_1, \mu_2,\cdots, \mu_m\}$ 
and the positive eigenvalues $\{\mu_n: n \geq m+1\}$ of $A$:
$$ H=H^+ \oplus H^-$$
where $H^+$ is a closed linear subspace of $H$ and $H^-$ is a
finite-dimensional subspace. Denote by $p^+ :H \to H^+$ and $p^-: H \to H^-$ the 
corresponding
projections onto $H^+$ and $H^-$ respectively. Since $H^-$ is 
finite-dimensional, then $T_t | H^-$ is invertible for each $t \geq 0$. 
Therefore, we can set  $T_{-t}:= [T_t|H^-]^{-1}: H^- \to H^-$ for each $t \geq 0$. 

 Consider the following semilinear stochastic evolution equation on H:
$$
\align
du(t)&=[-Au(t)+F(u(t))]\, dt+B_0dW(t), \quad t \geq 0, \tag{4.1}\\
u(0)&=x\in H.
\endalign 
$$
  In the above equation, let $F: H \to H$ be a globally Lipschitz 
map with Lipschitz constant $L$:
$$ |F(v_1)-F(v_2)| \leq L |v_1-v_2|, \quad v_1,v_2 \in H.$$
Then (4.1) has a unique mild solution given by
$$
u(t,x)=T_tx+\int _{0}^t T_{t-s}F(u(s,x))ds+\int _0^t T_{t-s}B_0dW(s), \quad t  \geq 0
\tag{4.2}
$$
Furthermore, if $F: H \to H$ is $C^{k,\e}$, the mild solution of (4.2) generates a 
$C^{k,\e}$ perfect cocycle also denoted by $u: \R^+ \times H \times \O \to H$.

Suppose that $F: H \to H$ is globally bounded, and its Lipschitz constant 
$L$ satisfies
$$ L[ \mu_{m+1}^{-1} - \mu_m^{-1}] < 1. \tag{4.3}$$
Note that the above condition is automatically satisfied in the affine linear 
case $F \equiv 0$. 

The next proposition is key to the existence and uniqueness of a stationary random 
point 
for the
cocycle $(u,\theta)$ in the sense of Definition 1.1.

\proclaim{Proposition 4.1}

Assume the above conditions on $A, B_0, F$ together with (4.3). Then there is a 
unique $\Cal F$-measurable map $Y: \O \to H$ satisfying 
$$
\aligned
Y(\omega)=\int _{-\infty}^0T_{-s}p^+F(Y(\theta (s,\omega)))ds-\int _0^{\infty}
T_{-s}p^-F(Y(\theta (s,\omega)))ds&\\
+(\o)\int _{-\infty}^0T_{-s}p^+B_0dW(s)-(\o)\int _0^{\infty}T_{-s}p^-B_0dW(s)
\endaligned
\tag{4.4}
$$
for all $\o \in \O$.
\endproclaim

\demo{Proof}

We use a contraction mapping argument to show that the integral equation 
(4.4) has an $\Cal F$-measurable solution $Y: \O \to H$. 

Define the $\Cal F$-measurable map $Y_1:  \O \to H$ by 
$$
Y_1(\o):=(\o)\int _{-\infty}^0T_{-s}p^+B_0\, dW(s)-(\o)\int 
_0^{\infty}T_{-s}p^-B_0\, 
dW(s), \quad \o \in \O.
$$
Denote by $B(\O,H)$ the Banach 
space of all (surely) bounded $\Cal F$-measurable maps $Z: \O \to H$ given the 
supremum norm $ \|Z\|_{\i}:= \displaystyle \sup_{\o \in \O} |Z(\o)|$.
Define the map $ \Cal M :B(\O,H) \to L^0 (\O,H)$ by 
$$
\align
\Cal M (Z)(\o)&:= \int _{-\infty}^0 T_{-s}p^+F(Z(\theta (s,\omega))+ Y_1(\theta 
(s,\omega)))\,ds\\
& \qquad \qquad \qquad -\int _0^{\infty}T_{-s}p^-F(Z(\theta 
(s,\omega))+Y_1(\theta 
(s,\omega)))\,ds
\endalign  
$$
for all $Z \in B(\O,H)$ and all $\o \in \O$. 

Note first that $\Cal M$ maps $B(\O,H)$ 
into itself. To see this let $Z \in B(\O,H)$ and $\o \in \O$. Then 
$$ 
\align
|\Cal M (Z)(\o)| &\leq \|F\|_{\i} \biggl [\int _{-\infty}^0 \|T_{-s}p^+\| \, ds 
+ \int _0^{\infty} \|T_{-s}p^-\|\,ds \biggr ]\\
&\leq \|F\|_{\i}\biggl [\int _{-\infty}^0 e^{s\mu_{m+1}} \, ds 
 + \int _0^{\infty}e^{s\mu_m }\,ds \biggr ] \\
&\leq \|F\|_{\i}[\mu_{m+1}^{-1}-\mu_m^{-1}] < \i 
\endalign
$$
where $\|F\|_{\i}:= \displaystyle \sup_{v \in H} |F(v)|$. Hence $\Cal M (Z) \in 
B(\O,H)$ for all $Z \in B(\O,H)$.  

Secondly, $\Cal M$ is a contraction.
To prove this, take any $Z_1, Z_2 \in 
B(\O,H)$ and $\o \in \O$. Then from the definition of $\Cal M$, we get
$$
\align
|\Cal M (Z_1)(\o)-\Cal M (Z_2)(\o)| &\leq L \int _{-\infty}^0 \|T_{-s}p^+\|\cdot 
|Z_1(\theta (s,\omega)) -Z_2(\theta (s,\omega))| \, ds \\
&+ L \int _0^{\infty} \|T_{-s}p^-\| \cdot |Z_1(\theta (s,\omega)) -Z_2(\theta 
(s,\omega))|   \,ds\\
&\leq L \|Z_1 -Z_{2}\|_{\i}\biggl [\int _{-\infty}^0 \|T_{-s}p^+\| \, ds + 
\int _0^{\infty} \|T_{-s}p^-\|\,ds \biggr ]\\
&\leq L \|Z_1 -Z_{2}\|_{\i} \biggl [\int _{-\infty}^0 e^{s\mu_{m+1}} \, ds 
 + \int _0^{\infty}e^{s\mu_m }\,ds \biggr ] \\
&=L [\mu_{m+1}^{-1}-\mu_m^{-1}]\|Z_1 -Z_{2}\|_{\i}\\
&=\mu \|Z_1 -Z_{2}\|_{\i}
\endalign
$$
where $\mu:=L [\mu_{m+1}^{-1}-\mu_m^{-1}] < 1$. This proves that 
$ \Cal M :B(\O,H) \to B (\O,H)$ is a contraction, and hence has a unique fixed 
point $Z_0 \in B (\O,H)$. That is 
$$
\align
Z_0 (\o)&:= \int _{-\infty}^0 T_{-s}p^+F(Z_0 (\theta (s,\omega))+ Y_1(\theta 
(s,\omega)))\,ds\\
& \qquad \qquad \qquad \qquad -\int _0^{\infty}T_{-s}p^-F(Z_0 (\theta 
(s,\omega))+Y_1(\theta (s,\omega)))\,ds
\endalign  
$$
for all $\o \in \O$. Now define $Y: \O \to H$ by 
$$ Y(\o):= Z_0 (\o) + Y_1 (\o), \quad \o \in \O.$$
It is easy to check that $Y$ satisfies the identity (4.4). 

Since $Z_0$ is uniquely determined, then so is $Y$. 
\qed

\enddemo

The following proposition gives existence and uniqueness 
of a stationary point for the see (4.1).

\proclaim{Proposition 4.2}

 Assume all the conditions on $A, B_0, F$ stated in Proposition 4.1. Suppose 
that $F$ is globally bounded, globally Lipschitz  and satisfies condition (4.3). 
Then the semilinear see (4.1) has a unique 
stationary point $Y: \O \to H$, i.e. 
$u(t,Y(\omega),\omega)=Y(\theta (t,\omega))$ for all $t \geq 0$ and $\o \in \O$. 
Furthermore,
$Y \in L^p (\O,H)$ for all $p \geq 1$.
\endproclaim


\demo{Proof}

By hypotheses and Proposition 4.1, the integral equation (4.4)  has a unique $\Cal 
F$-measurable solution $Y: \O \to H$. Let $t \geq 0$. Using (4.4), it follows that 
$$
\allowdisplaybreaks
\align
 Y(\theta (t,\omega))&=\int 
_{-\infty}^0T_{-s}p^+F(Y(\theta (t+s,\omega)))ds-\int _0^{\infty}
T_{-s}p^-F(Y(\theta (t+s,\omega)))\,ds\\
&+(\o) \int _{-\infty}^0T_{-s}p^+B_0\, dW(s+t)-(\o)\int 
_0^{\infty}T_{-s}p^-B_0\, dW(s+t)\\
=&\int_{-\infty}^tT_{t-s}p^+F(Y(\theta (s,\omega)))ds-\int _t^{\infty}
T_{t-s}p^-F(Y(\theta (s,\omega))) \,ds\\
&+(\o)\int _{-\infty}^tT_{t-s}p^+B_0\, dW(s)-(\o)\int _t^{\infty}T_{t-s}p^-B_0\, 
dW(s)\\
=&T_t \biggl [\int_{-\infty}^0T_{-s}p^+F(Y(\theta (s,\omega)))ds-\int 
_0^{\infty}
T_{-s}p^-F(Y(\theta (s,\omega)))\,ds\\
&+(\o)\int _{-\infty}^0T_{-s}p^+B_0\, dW(s)-(\o)\int _0^{\infty}T_{-s}p^-B_0\, 
dW(s)\biggr 
]\\
&+\int _{0}^t T_{t-s}p^+F(Y(\theta (s,\omega))ds+\int _0^t
T_{t-s}p^-F(Y(\theta (s,\omega)))\,ds\\
&+(\o)\int _0^t T_{t-s}p^+B_0\, dW(s)+(\o)\int _0^tT_{t-s}p^-B_0\, dW(s)\\
=&T_tY(\omega)+\int _{0}^t 
T_{t-s}F(Y(\theta (s,\omega)))ds + (\o)\int _0^t T_{t-s}B_0\, dW(s).
\endalign
$$
This gives
$$
Y(\theta (t,\omega))=T_tY(\omega)+\int _{0}^t 
T_{t-s}F(Y(\theta (s,\omega)))ds+ (\o)\int _0^t T_{t-s}B_0\, dW(s)
$$
for all $t \geq 0$. Therefore, $Y(\theta (t,\o)), t \geq 0, \o \in \O,$ is a 
stationary solution of (4.2) (with $x=Y(\o)$). Since $u(t,Y(\omega),\omega), t 
\geq 0, \o \in \O$, is also a solution of (4.2), then by 
uniqueness of the solution to (4.2), we must have
$$
u(t,Y(\omega),\omega)=Y(\theta (t,\omega))   
$$
for all $t \geq 0$ and all $\o \in \O$. Hence $Y$ is a stationary point for the 
see (4.1).  

The stationary point for (4.1) is unique (within the class of  $\Cal F$-measurable 
maps 
$\O \to H$). To see this,
it is sufficient to observe that the above computation shows that every stationary 
point 
of (4.1) is a solution of the integral
equation (4.4). Uniqueness of the stationary solution then follows from Proposition 
4.1. 

In view of the proof of Proposition 4.1, the last assertion of Proposition 4.2 follows 
from the fact that 
$Y_1 \in L^p (\O,H)$ for all $p \geq 1$ and $Z_0 \in L^\i (\O,H)$ .
\qed.
\enddemo

 \medskip

The existence of local stable and unstable manifolds near a stationary point of the 
affine stochastic evolution 
equation (4.1) follows from a straightforward modification of the proof of Theorem 4.1 
in 
the next section.

\bigskip
\noindent
{\it (b) Semilinear stochastic evolution equations: Linear noise}

Here we recall the setting and hypotheses leading to Theorem 2.6 in [M-Z-Z.1]. 

We will prove the existence of local stable and unstable manifolds  for semiflows 
generated by mild 
solutions of semilinear stochastic evolution equations of the form:
$$
\left. \aligned
du(t)&=-Au(t) dt+F(u(t))dt+ Bu(t)\,dW(t), \quad t > 0,\\
u(0)&=x \in H.
\endaligned \right \} \tag{4.5}
$$
%
 
In the above equation  $A: D(A) \subset H \to H$ is a
closed linear operator on a separable real Hilbert space $H$. Assume 
that $A$ has a complete orthonormal system of eigenvectors  $\{ e_n: n\geq 1\}$ with 
corresponding positive
eigenvalues $\{ \mu_n, n\geq 1\}$;  i.e.,
$Ae_n=\mu_n e_n, \,\, n \geq 1.$
Suppose $-A$ generates  a strongly continuous semigroup
of bounded linear operators $T_t: H \to H, \, t\geq 0$. Let $E$ be a separable 
real Hilbert space.
Suppose $W(t), t\geq 0,$ is $E$-valued cylindrical
Brownian motion defined on the canonical filtered Wiener space
$(\O, \F, (\F_t)_{t\geq 0},P)$
and with  a separable covariance Hilbert space $K$, where $K \subset E$ is 
a Hilbert-Schmidt embedding. That is, $\O$ is the space
of all continuous paths $\o: \R \to E$ such that $\o (0)=0$ with the 
compact open
topology, $\F$ is its Borel $\sigma$-field,  $\F_t$ is the sub-$\sigma$-field
generated by all evaluations $\O \ni \o \mapsto \o(u) \in E, u \leq t$, and $P$
is Wiener measure on $\O$. The Brownian motion is given by
$$ W(t,\o):=\o(t), \quad \o \in \O,\, t \in \R,$$
and may be represented by
$$
W(t)=\sum_{k=1}^{\i} W^k(t) f_k, \quad t \in \R,
$$
  where $\{f_k: k \geq 1\}$ is a complete orthonormal basis of $K$, and the
$W^k, k \geq 1,$ are standard independent one-dimensional Wiener 
processes ([D-Z.1],
Chapter 4).

Suppose $B: H \to L_2(K,H)$ is  a bounded linear operator. The stochastic
integral in (4.5) is defined in the sense of ([D-Z.1], Chapter 4). 

Assume the hypotheses of Theorem 2.4 in [M-Z-Z.1].

We will denote by
$\theta: \R \times \O \to \O$ the standard $P$-preserving ergodic 
Wiener shift on
$\O$:
$$\theta (t,\o)(s):=\o (t+s)-\o(t), \quad t,s \in \R.$$
 Let $L(H)$ be the Banach space of all bounded linear operators $H
\to H$ given the
uniform operator norm $\|\cdot\|$.  Denote by $L_2(H) \subset L(H)$ the Hilbert
space  of all
Hilbert-Schmidt operators $S: H \to H$. 

Suppose $F:H \to H$ is a (Fr\'echet) $C^{k,\e}$ ($k \geq 1, \e \in (0,1]$) non-linear 
map satisfying the following Lipschitz
and linear growth hypotheses:
$$
\left. \aligned
|F(v)| &\leq C(1+|v|), \quad v \in H\\
|F(v_1)-F(v_2)| &\leq L_n |v_1-v_2|, \quad v_i \in H,
|v_i| \leq n, i=1,2,
\endaligned \right \} \tag{4.6}
$$
for some positive constants $C,L_n, n \geq 1 $.

The mild solutions of the see (4.5) generate a $C^{k,\e}$ ($k \geq 1, \e \in (0,1]$) 
perfect cocycle 
$(U, \theta)$ on $H$, satisfying all the assertions of Theorem 2.6 of [M-Z-Z.1].

Under the above conditions, one gets the following stable manifold theorem 
for hyperbolic stationary trajectories of the see (4.5). 

 \proclaim{Theorem 4.1}

Assume the above hypotheses on the coefficients of the see (4.5). Assume that the 
stochastic semiflow 
$U: \R^+ \times H \times \O \to H$ generated by mild solutions of (4.5) has a 
hyperbolic 
stationary point $Y:\O \to H$ 
such that $E \log^+ |Y| < \i$.  Then $(U,\theta)$ has a perfect family of $C^{k,\e}$ 
local 
stable and unstable manifolds 
satisfying all the assertions of Theorem 2.1. 
\endproclaim

\demo{Proof}

One first checks the estimate
$$ \int_{\O} \log^+ \sup_{0\leq t_1,t_2 \leq a} \|U(t_2, Y(\theta
(t_1,\o))+(\cdot),\theta (t_1,\o))\|_{k,\e} \, dP(\o) < \infty \tag{4.7}$$
for any fixed $0 <  \rho, a < \infty, k \geq 1$  and $\e \in (0,1]$. This 
estimate 
follows from the integrability condition on $Y$ and assertion (vi) 
of Theorem 2.6 in [M-Z-Z.1]. The conclusion of Theorem 4.1 now follows 
immediately from Theorem 2.1. \qed  
\enddemo

\bigskip
\noindent
{\it (c) Semilinear parabolic spde's: Lipschitz nonlinearity}

Consider the Laplacian 
$$
\Delta:={1\over 2}\sum \limits _{i,j=1}^d  {\partial ^2\over
\partial \xi_i^2 }  \tag{4.8}
$$
defined on a smooth bounded domain $\Cal D$ in $\R^d$, with a smooth boundary 
$\partial 
\Cal D$
with zero Dirichlet boundary conditions.
Assume  that $f: \R \to \R$ is a $C^{\i}_b$ function and let $d\xi$ be 
Lebesgue measure on $\R^d$. 
Let $W_n, n \geq 1,$ be
independent one-dimensional standard Brownian motions with 
$W_n(0)=0$ defined on the canonical filtered Wiener space $(\O,\F, P, (\F_t)_{t \in 
\R})$. Let $\theta$ denote the Brownian
shift on $\O:= C(\R, \R^\i)$. 
Recall that the Sobolev 
space  $H_0^{k} (\Cal D)$ is
 the completion of $C_0^{\infty}(\Cal D,\R)$ under the Sobolev norm
$$
\vert \vert u\vert \vert_{H_0^{k} (\Cal D)}^2 :=\sum_{|\alpha |\leq k}\int_{\Cal D} 
|D^{\alpha}u(\xi)|^2 \, d\xi.
$$
Suppose further that $\sigma_n \in H_0^{s} (\Cal D), n \geq 1$, and  the series 
$\displaystyle \sum_{n=1}^\i \sigma_n$
converges absolutely 
in $H_0^{s} (\Cal D)$ where $s > k + \displaystyle \frac{d}{2} > d$.

By Theorem 3.5 ([M-Z-Z.1]), weak solutions of the initial-value problem:  
$$
\left. \aligned
du(t)&= \frac{1}{2} \Delta u(t) dt+f(u(t))dt+\sum_{n=1}^\i \sigma_n(\xi)u(t)\,dW 
_n(t), \quad 
t > 0\\
u(0)&=\psi \in  H_0^{k} (\Cal D) 
\endaligned \right \} \tag{4.9}
$$
give a perfect smooth cocycle $(U,\theta)$ on the Sobolev space $H_0^{k}(\Cal D)$  
which satisfies all the assertions of
Theorem 3.5 in [M.Z.Z.1].   Applying Theorem 2.1, we get the following stable manifold 
theorem for the spde (4.9):

   %
   %

\proclaim{Theorem 4.2}

Assume the above hypotheses on the coefficients of the spde (4.9). Assume that the 
stochastic semiflow $U: \R^+ \times
H_0^{k} (\Cal D)  \times \O \to   H_0^{k} (\Cal D)$    generated by weak solutions of 
(4.9) 
has a hyperbolic stationary point $Y:\O
\to H_0^{k} (\Cal D)$   such that $E \log^+ \|Y\|_{H_0^{k} }  < \i$.  Then 
$(U,\theta)$ 
has a perfect family of $C^\i$ local stable
and unstable manifolds in $H_0^{k} (\Cal D)$  satisfying all the assertions of Theorem 
2.1. 
\endproclaim

\bigskip
\noindent
{\it (d) Stochastic reaction diffusion equations: dissipative nonlinearity}


In section 4 (a) of [M.Z.Z.1], we constructed a $C^1$ stochastic semiflow on the 
Hilbert 
space $H:=L^2(\Cal D)$
 for a stochastic reaction-diffusion equation 
$$
du=\nu \Delta u\, dt+u(1-|u|^{\alpha})\,dt+\sum_{i=1}^\i \sigma_i(\xi)u(t)\,dW _i(t),
\tag4.10
$$
defined on a bounded domain $\Cal D \subset \R^d$ with a smooth boundary $\partial 
\Cal 
D$.  In (4.10), 
the Laplacian on $\Cal D$ is
denoted by $\Delta$, and we impose  Dirichlet boundary conditions on $\partial \Cal 
D$. 
The  $W_i, i \geq 1,$  are independent one-dimensional standard Brownian 
motions and $\sum \limits
_{i=1}^{\infty}\sigma_i$ is absolutely convergent in $H^s(\Cal D)$, for $s>2+{d\over 
2}$. The dissipative term
yields the existence of a unique stationary solution of (4.10) under a suitable choice 
of 
the parameter $\nu$ ([D-Z.2]). 

In view of the estimates in Theorem 4.1 ([M.Z.Z.1]) and Theorem 2.1, one gets the 
following:

\proclaim{Theorem 4.3}

Assume the above hypotheses on the coefficients of the spde (4.10). Let 
$\alpha < {4 \over 
d}$. Assume that the stochastic semiflow $U: \R^+ \times
L^2(\Cal D) \times \O \to  L^2(\Cal D) $  generated by mild solutions of (4.10) has a 
hyperbolic stationary point 
$Y:\O \to L^2(\Cal D)$   such that $E \log^+ \|Y\|_{L^2}  < \i$.  Then $(U,\theta)$ 
has 
a 
perfect family of $C^1$ local stable and
unstable manifolds in   $L^2(\Cal D)$   satisfying all the assertions of Theorem 2.1. 
\endproclaim



\remark{Remarks}

\item{(i)} The results in Sections  (c) and (d) hold if the Euclidean domain $\Cal D$ 
is 
replaced by a compact smooth
$d$-dimensional Riemannian manifold $M$ 
(possibly with a smooth boundary $\partial M$).

\item{(ii)} We conjecture that Theorem 4.3 still holds (but with {\it Lipschitz} 
stable/unstable 
manifolds) if the dissipative term $u(1-|u|^{\alpha})$ is replaced by a more general 
one 
of the
form $F(u):=f\circ u$, where $f: \R \to \R$ is a $C^1$ function satisfying the 
following 
classical estimates:
$$
\align
-c_1- \alpha_1 |x|^p \leq &f(x)x \leq c_1 -\alpha_2 |x|^p,  \quad f'(x) \leq  c_2,
\endalign
$$
for all $x \in \R$, with $c_1, c_2, \alpha_1, \alpha_2$ positive constants, and $p$ 
any 
integer greater than $2$.

\item{(iii)} Is it true that the stochastic flow and the local stable/unstable 
manifolds in Theorem 4.3 are of class $C^2$?

\endremark

\bigskip
\noindent
{\it (e) Stochastic Burgers equation: additive noise}

The existence of a $C^1$ stochastic semiflow on $L^2([0,1])$ for Burgers equation
$$
du+u{\partial u \over \partial \xi }\,dt=\nu \Delta udt+dW(t), \quad t > 0, \, \, \nu 
> 0,
\tag{4.11}
$$
was established in Part I of this work, where $W(t), \, t > 0,$ is an infinite 
dimensional 
Brownian motion on $L^2[0,1]$. See [M-Z-Z.1], Theorem 4.3. 

Under extra spatial smoothness hypotheses on the noise, viz. $W(t,\cdot) \in 
C^3([0,1])$, Burgers equation (4.11) admits a unique stationary point ([Si]).   More 
generally, with our weaker condition on the noise $W$ ([M.Z.Z.1], Section 4 (b)), we 
stipulate that
equation (4.11) has a hyperbolic stationary point. In this case, we get the 
following result:
 
 \proclaim{Theorem 4.4}

Assume the  hypotheses of Theorem 4.3  of [M.Z.Z.1] on the coefficients of Burgers 
spde (4.11). 
Assume that the stochastic semiflow $U: \R^+ \times L^2 ([0,1]) \times \O \to L^2 
([0,1]) $  
generated by mild solutions
of (4.11) has a hyperbolic stationary point 
$Y:\O \to L^2 ([0,1]) $  such that $E \log^+ \|Y\|_{L^2}  < \i$.  Then $(U,\theta)$ 
has 
a perfect family of $C^1$ local stable
and unstable manifolds in   $L^2([0,1])$  satisfying all the assertions of Theorem 
2.1.  
\endproclaim

\medskip
 
Note that hyperbolicity of the stationary point in Theorem 4.4 is in the sense of 
Definition (1.3).  
Theorem 1.1 of this article and Theorem 4.3 of ([M.Z.Z.1]) imply  that the
Lyapunov spectrum for the linearization of (4.11) exists and is discrete for any 
viscosity $\nu > 0$. 
When $W$ is  $C^3$ in the space variable, it known that for any
$C^2$ initial condition, the solution $u(t)$ of (4.11)  converges to the stationary 
solution for any positive viscosity $\nu > 0$ ([Si]). It is therefore easy to see 
that the stable manifold is the whole of $L^2[0,1]$. 
  
The case of sufficiently large viscosity and rough noise $W(t) \in L^2([0,1])$ 
is currently being studied ([L-Z]). This work 
shows  that (4.11) admits a unique globally exponentially 
stable stationary point in this case.  So in this (somewhat non-generic) case, the 
unstable 
manifold consists of the single random point $Y(\o) \in L^2([0,1])$, and the  
non-linear cocycle $U$ will approach $Y(\o)$ with exponential speed less than or equal 
to the top Lyapunov exponent $\lambda_1$ of the linearized Burgers equation. 

We conjecture that the assertions in the above paragraph still hold for any viscosity 
$\nu > 0$
(cf. [D-Z.2], Theorem 14.4.4). Further analysis of the  Lyapunov spectrum for (4.11) 
(in the cases of small and zero 
viscosity $\nu$) is postponed to a future project.

\bigskip

\centerline{\bf Acknowledgments}
\medskip

The authors would like to thank Prof. B. $\emptyset$ksendal for inviting them
to Oslo in the summer 2000 where the project
was started. They are also grateful to Prof. K.D. Elworthy for
inviting them to Warwick
on numerous occasions especially during the Warwick SPDE's Symposium
2000/2001 so that they may have
opportunities to meet; to Prof. A. Truman for inviting them to the
International
Workshop of Probabilistic Methods in Fluids at Swansea in April 2002
where  preliminary versions of the results
were announced; and to the organizers of ICM 2002 and First
Sino-German Stochastic Analysis Conference
in Beijing in August 2002 where the results were presented.
We would also like to thank Z. Brzezniak, F. Flandoli, Y. Liu, J. Robinson, J. 
Zabczyk,  and 
several other colleagues 
for useful discussions.

\bigskip

\newpage

\centerline{\bf REFERENCES}

\bigskip
\baselineskip=14truept
\parskip=3truept

\medskip







\item{[Co]\quad} Cohn, D\. L\., {\it Measure Theory\/}, Birkh{\"a}user (1980).
 \medskip


 \item{[C-V]\quad} Castaing, C\., and Valadier, M\., {\it Convex Analysis and 
Measurable 
Multifunctions}, Lecture Notes in Mathematics, vol. 580, Springer-Verlag, 
Berlin-Heidelberg-New York (1977).

\medskip

\item{[D-Z.1]\quad} Da Prato, G., and Zabczyk, J., {\it 
Stochastic Equations in Infinite Dimensions}, Cambridge University 
Press (1992).

\medskip

\item{[D-Z.2]\quad} Da Prato, G., and Zabczyk, J., {\it Ergodicity for Infinite 
Dimensional Systems}, 
Cambridge University Press (1996).

\medskip
\item{[F-S]\quad}Flandoli, F\., and Schauml\"offel, K\.-U\., A multiplicative ergodic 
theorem with applications to a first order
stochastic hyperbolic equation in a bounded domain, {\it  Stochastics and Stochastics 
Reports\/},  34  (1991),  no. 3-4,
241--255.

 \medskip

\item{[L-Z]\quad} Liu, Y. and Zhao, H.Z., The stationary solutions of some SPDEs, in 
preparation.
 
\medskip
 
 \item{[Mo.1]\quad}  Mohammed, S\.-E\. A\., The Lyapunov spectrum and 
stable manifolds for stochastic linear delay equations, 
{\it Stochastics and Stochastic Reports\/}, Vol. 29 (1990), 89-131.

\medskip

 


\item{[M-S.1]\quad}  Mohammed, S\.-E\. A\., and  Scheutzow, M\. K\. R\.,
Lyapunov exponents of linear stochastic functional differential equations driven by  
semimartingales, Part I:  The multiplicative ergodic theory, 
 {\it Ann. Inst. Henri Poincar{\'e}, Probabilit{\'e}s et Statistiques,\/}
 Vol. 32, 1, (1996), 69-105. pp. 43.

\medskip
 \item{[M-S.2]\quad}  Mohammed, S\.-E\. A\., and  Scheutzow, M\. K\. R\.,
 Spatial 
estimates for stochastic flows in Euclidean space,   
 {\it The Annals of Probability}, 26, 1, (1998), 56-77. 

\medskip

\item{[M-S.3]\quad}  Mohammed, S\.-E\. A\., and  Scheutzow, M\. K\. R\.,
 The stable manifold theorem for stochastic differential equations, 
{\it The Annals of Probability\/}, Vol. 27, No. 2, (1999),  615-652.     

\medskip

\item{[M-S.4]\quad}  Mohammed, S\.-E\. A\., and  Scheutzow, M\. K\. R\.,
 The stable manifold theorem for non-linear stochastic systems with memory.
Part I: Existence of the semiflow, Part II:  The local stable manifold theorem. 
 (preprints).
 
 \medskip
  
\item{[M-Z-Z.1]\quad}  Mohammed, S\.-E\. A\., Zhang, T\. S\. and Zhao, H\.,
The stable manifold theorem for semilinear stochastic evolution equations and 
stochastic partial differential equations I: Flows and stationary solutions (preprint) 
(2002).
\medskip

\item{[O]\quad} Oseledec, V\. I\.,  A multiplicative ergodic theorem.
Lyapunov characteristic numbers for dynamical systems, {\it Trudy Moskov.
Mat. Ob\v s\v c.\/} 19 (1968), 179-210. English transl. {\it Trans. Moscow
Math. Soc.\/} 19 (1968), 197-221.  
\medskip

\item{[Ru.1]\quad} Ruelle, D., Ergodic theory of differentiable dynamical
systems,
{\it Publ. Math. 
Inst. Hautes Etud. Sci.} (1979), 275-306.

\medskip
\item{[Ru.2]\quad} Ruelle, D\., Characteristic exponents and invariant
manifolds in
Hilbert space, {\it Annals of Mathematics 115\/} (1982), 243--290.

\medskip



\item{[Si]\quad} Sinai, Ya. G., Burgers system driven by a periodic stochastic flow,
In: {\it It\^o's stochastic calculus and probability theory},  Springer, Tokyo (1996), 
347--353.

\bigskip

\noindent

\bigskip

\noindent
Salah-Eldin A. Mohammed \newline
Department of Mathematics,
\newline
Southern Illinois University at Carbondale,
\newline
Carbondale, Illinois 62901.\newline
Email: salah\@sfde.math.siu.edu \newline
Web page: http://sfde.math.siu.edu

\bigskip
\medskip
\noindent
Tusheng Zhang \newline
Department of Mathematics \newline
     University of Manchester,\newline
Oxford Road, Manchester M13 9PL\\ \newline
UK.
     \newline
Email: tzhang\@math.man.ac.uk

\bigskip
\medskip
\noindent
Huaizhong Zhao \newline
Department of Mathematical Sciences \newline
Loughborough University,\newline
LE11 3TU,\newline
UK.
     \newline
Email: H.Zhao\@lboro.ac.uk

\end

\enddocument